\documentclass[11pt]{article}
\usepackage{amsfonts, amsmath,amssymb}
\usepackage{hyperref}
\usepackage{graphicx}

\newcommand{\detail}[1]{\par\noi{\bf [Proof detail\ }{#1}
\hfill{\bf ]}\par\noi\hspace{-4pt}}
\renewcommand{\detail}[1]{}

\newcommand{\file}{wisk/siva/surv/surf.tex\quad}
\renewcommand{\file}{}

\newcommand{\plaat}[2]{{#2}}
\newcommand{\dis}{\displaystyle}

\newcommand{\noi}{\noindent}

\newcommand{\halmos}{\rule{1ex}{1.4ex}}
\def \qed {\nopagebreak{\hspace*{\fill}$\halmos$\medskip}}
\newcommand{\med}{\medskip}

\newtheorem{theorem}{Theorem}
\newtheorem{proposition}[theorem]{Proposition}
\newtheorem{corollary}[theorem]{Corollary}
\newtheorem{conjecture}[theorem]{Conjecture}
\newtheorem{lemma}[theorem]{Lemma}

\newtheorem{remark}[theorem]{Remark}

\newcommand{\bt}{\begin{theorem}}
\newcommand{\et}{\end{theorem}}
\newcommand{\bl}{\begin{lemma}}
\newcommand{\el}{\end{lemma}}
\newcommand{\bp}{\begin{proposition}}
\newcommand{\ep}{\end{proposition}}
\newcommand{\bcor}{\begin{corollary}}
\newcommand{\ecor}{\end{corollary}}
\newcommand{\br}{\begin{remark}\rm}
\newcommand{\er}{\end{remark}}
\newcommand{\bcon}{\begin{conjecture}}
\newcommand{\econ}{\end{conjecture}}

\newcommand{\be}{\begin{equation}}
\newcommand{\ee}{\end{equation}}
\newcommand{\ba}{\begin{array}}
\newcommand{\ea}{\end{array}}
\newcommand{\bc}{\be\begin{array}{r@{\,}c@{\,}l}}
\newcommand{\ec}{\end{array}\ee}

\newcommand{\al}{\alpha}
\newcommand{\bet}{\beta}
\newcommand{\ga}{\gamma}

\newcommand{\de}{\delta}
\newcommand{\De}{\Delta}
\newcommand{\eps}{\varepsilon}
\newcommand{\la}{\lambda}
\newcommand{\La}{\Lambda}

\newcommand{\tet}{\theta}

\newcommand{\om}{\Omega}

\newcommand{\Fi}{{\cal F}}

\newcommand{\Hi}{{\cal H}}

\newcommand{\R}{{\mathbb R}}

\newcommand{\Z}{{\mathbb Z}}

\renewcommand{\P}{{\mathbb P}}
\newcommand{\E}{{\mathbb E}}

\newcommand{\desd}{\ensuremath{\Leftrightarrow}}

\newcommand{\sub}{\subset}
\newcommand{\beh}{\backslash}

\newcommand{\ti}{\tilde}
\newcommand{\dgg}{\dagger}
\newcommand{\ov}{\overline}
\newcommand{\un}{\underline}
\newcommand{\subb}[2]{_{\ba{c}\scriptstyle{#1}\\[-.15cm]\scriptstyle{#2}\ea}}

\newcommand{\ffrac}[2]{{\textstyle\frac{{#1}}{{#2}}}}
\newcommand{\dif}[1]{\ffrac{\partial}{\partial{#1}}}
\newcommand{\diff}[1]{\ffrac{\partial^2}{{\partial{#1}}^2}}

\newcommand{\half}{{[0,\infty)}}
\newcommand{\expo}{\mbox{\large\it e}}
\newcommand{\ex}[1]{\expo^{\,\textstyle{#1}}}

\newcommand{\comb}{\,\&\,}

\begin{document}

%numbering formulas within sections
\makeatletter\@addtoreset{equation}{section}
\makeatother\def\theequation{\thesection.\arabic{equation}} 

%alternative layout for enumerate lists.
\renewcommand{\labelenumi}{{(\roman{enumi})}}

\title{\vspace{-3cm}Survival of contact processes on the hierarchical group}
\author{
Siva R.~Athreya\vspace{6pt}\\
{\small Indian Statistical Institute}\\
{\small 8th Mile, Mysore Road}\\
{\small RVCE Post}\\
{\small Bangalore 560 059}\\
{\small India}\\
{\small e-mail: athreya@isibang.ac.in}\vspace{8pt}
\and Jan M.~Swart\vspace{6pt}\\
{\small Institute of Information}\\ 
{\small Theory and Automation}\\
{\small of the ASCR (\' UTIA)}\\
{\small Pod vod\'arenskou v\v e\v z\' i 4}\\
{\small 18208 Praha 8}\\
{\small Czech Republic}\\
{\small e-mail: swart@utia.cas.cz}\vspace{4pt}}
\date{{\small\file} March 2, 2009}

\maketitle\vspace{-30pt}

\begin{abstract}\noi
We consider contact processes on the hierarchical group, where sites infect
other sites at a rate depending on their hierarchical distance, and sites
become healthy with a constant recovery rate. If the infection rates decay too
fast as a function of the hierarchical distance, then we show that the
critical recovery rate is zero. On the other hand, we derive sufficient
conditions on the speed of decay of the infection rates for the process to
exhibit a nontrivial phase transition between extinction and survival. For our
sufficient conditions, we use a coupling argument that compares contact
processes on the hierarchical group with freedom two with contact processes on
a renormalized lattice. An interesting novelty in this renormalization
argument is the use of a result due to Rogers and Pitman on Markov functionals.
\end{abstract}

\noi
%82C22 Interacting particle systems
%60K35 Interacting random processes; statistical mechanics type models,
%      percolation theory
%82C28 Dynamic renormalization group methods
%
%82C26 Dynamic and nonequilibrium phase transitions (general)
{\it Acknowledgements.} Work sponsored by GA\v CR grant 201/06/1323 and CSIR
Grant 25(0164) /08/EMRII. Much of this work was carried out during a visit of
J.M.~Swart to the Indian Statistical Institute, Bangalore centre and Kolkata
centre, during Nov.\ and Dec.\ 2007.

{\small\setlength{\parskip}{-2pt}\tableofcontents}
\newpage

\section{Introduction}

\subsection{Main result}

Let $\La$ be a finite or countably infinite set, called lattice, let
$(a(i,j))_{i,j\in\La,\ i\neq j}$ be nonnegative constants, and $\de\geq
0$. Then the {\em contact process} on $\La$ with {\em infection rates}
$a(i,j)$ and {\em recovery rate} $\de$ is the $\{0,1\}^\La$-valued Markov
process $X=(X_t)_{t\geq 0}$ with the following description. If $X_t(i)=0$
(resp.\ $X_t(i)=1$), then we say that the site $i\in\La$ is {\em healthy}
(resp.\ {\em infected}) at time $t\geq 0$. An infected site $i$ infects a
healthy site $j$ with rate $a(i,j)\geq 0$, and infected sites become healthy
with rate $\de\geq 0$. It can be shown (see \cite[Prop.~I.3.2]{Lig85}) that
$X$ is well-defined provided the infection rates are summable, in the sense
that
\be\label{summ}
|a|:=\sup_{j\in\La}\sum_{i\in\La,\ i\neq j}a(i,j)<\infty.
\ee
Usually, it is convenient to assume also that $|a^\dgg|<\infty$, where
$|a^\dgg|$ is defined as in (\ref{summ}), but for the reversed infection rates
$a^\dgg(i,j):=a(j,i)$.

A contact process may be used to model the spread of an infection in a
spatially ordered population; see \cite{Lig99} as a general reference. A basic
feature of the contact process is that it exhibits a phase transition between
survival and extinction. Let $0\in\La$ be some fixed site, called {\em
  origin}. We say that a contact process on a lattice $\La$ with given
infection rates $(a(i,j))_{i,j\in\La,\ i\neq j}$ and recovery rate $\de\geq 0$
{\em survives} if there is a positive probability that the process started
with only the origin infected never recovers completely, i.e., if
\be\label{Pde}
\P^{\de_0}[X_t\neq\un 0\ \forall t\geq 0]>0,
\ee
where $\de_i\in\{0,1\}^\La$ is defined as $\de_i(j):=1$ if $i=j$ and
$\de_i(j):=0$ otherwise, and $\un 0\in\{0,1\}^\La$ denotes the configuration
with only healthy sites. (In typical cases, e.g.\ when the infection rates are
irreducible in an appropriate sense or if the process has some
translation-invariant structure, this definition will not depend on the choice
of the origin $0$.)

For given infection rates, we let
\be\ba{r@{\,}l}
\dis\de_{\rm c}:=\sup\big\{\de\geq 0:&\dis\mbox{the contact process with
infection rates }\\
&\dis(a(i,j))_{i,j\in\La,\ i\neq j}
\mbox{ and recovery rate $\de$ survives}\big\}
\ec
denote the {\em critical recovery rate}. A simple monotone coupling argument
shows that $X$ survives for $\de<\de_{\rm c}$ and dies out for
$\de>\de_{\rm c}$.

By comparison with a subcritical branching process, it is not hard to show
that $\de_{\rm c}\leq|a^\dgg|$, where $|a^\dgg|$ is defined below
(\ref{summ}). For a large class of lattices, it is known that moreover
$\de_{\rm c}>0$. For example, this is the case for nearest-neighbor processes
on infinite graphs, where $a(i,j)$ equals some fixed constant $\la>0$ if $i$
and $j$ are connected by an edge and is zero otherwise, or if $\La$ is a
finitely generated, infinite group, and the infection rates are irreducible
and invariant under the left action of the group \cite[Lemma~4.18]{Swa07}. On
groups that are not finitely generated, the question whether $\de_{\rm c}>0$
becomes more subtle. Inspired by a question that came up in \cite{Swa08}, the
main aim of the present paper is to give sufficient conditions for $\de_{\rm
  c}>0$ (resp.\ $\de_{\rm c}=0$) when $\La$ is the hierarchical group.

By definition, the {\em hierarchical group with freedom $N$} is the set
\be\label{omdef}
\om_N:=\big\{i=(i_0,i_1,\ldots):
i_k\in\{0,\ldots,N-1\},\ i_k\neq 0\mbox{ for finitely many }k\big\},
\ee
equipped with componentwise addition modulo $N$. We set
\be\label{norm}
|i|:=\inf\{k\geq 0:i_m=0\ \forall m\geq k\},
\ee
and call $|i-j|$ the {\em hierarchical distance} between two elements
$i,j\in\om_N$. We will be interested in contact processes on $\om_N$ whose
infection rates $a(i,j)$ are a function of the hierarchical distance between
$i$ and $j$ only. Such infection rates may always be written as
\be\label{adef}
a(i,j):=\al_{|i-j|}N^{-|i-j|}\qquad(i,j\in\om_N,\ i\neq j)
\ee
where $(\al_k)_{k\geq 1}$ are nonnegative constants. The scaling with
$N^{-|i-j|}$ in (\ref{adef}) is chosen for calculational convenience.
It is easy to check that in order for the infection rates $a(i,j)$ to summable
in the sense of (\ref{summ}), we must assume that
$\sum_{k=1}^\infty\al_k<\infty$.

Here is our main result:
\bt{\bf((Non-) triviality of the critical recovery rate)}\label{T:main}
Let $N\geq 2$, let $(\al_k)_{k\geq 1}$ be nonnegative constants such that
$\sum_{k=1}^\infty\al_k<\infty$, and let $\de_{\rm c}$ be the critical death
rate of the contact process on $\om_N$ with infection rates as in
(\ref{adef}).\med

\noi
{\bf(a)} Assume that
\be\label{loglim}
\liminf_{k\to\infty}N^{-k}\log(\bet_k)=-\infty,\quad\mbox{where}
\quad\bet_k:=\sum_{n=k}^\infty\al_n\quad(k\geq 1).
\ee
Then $\de_{\rm c}=0$.\med

\noi
{\bf(b)} In case $N$ is a power of 2, assume that
\be\label{logsum}
\sum_{k=m}^\infty N^{-k}\log(\al_k)>-\infty\qquad\mbox{for some }m\geq 0.
\ee
Otherwise, assume that (\ref{logsum}) holds with $N$
replaced by some real $N'<N$. Then $\de_{\rm c}>0$.
\et
The special role played by powers of 2 in part~(b) is entirely due to our
methods of proof and has no real significance. In fact, we will carry out
most of our calculations for the case $N=2$ and then generalize to the
statement in part~(b) by a comparison argument. Note that if the
$\al_k$ have the  double exponential form
\be\label{examp}
\al_k=e^{-\tet^k}\qquad(k\geq 1),
\ee
then our results show that $\de_{\rm c}>0$ for $1<\tet<N$ and $\de_{\rm c}=0$
for $\tet>N$. There is a gap between the conditions (\ref{loglim}) and
(\ref{logsum}). We guess that (\ref{logsum}) is not necessary for $\de_{\rm
  c}>0$, since this condition is violated when infinitely many of the
$\al_k$'s are zero, while it seems unlikely that the latter should imply
$\de_{\rm c}=0$. We do not know if condition (\ref{loglim}) is sharp.

\detail{Indeed, if $\al_k$ is as in (\ref{examp}), and $\tet<N$, then we can
  choose $\tet<N'<N$ such that
\[
\sum_k(N')^{-k}\log(\al_k)=-\sum_k\big(\ffrac{\tet}{N'}\big)^k>-\infty,
\]
hence (\ref{logsum}) holds. (Note that we need $\tet>1$ in order to have
$\sum_k\al_k<\infty$.) Now
\[\ba{r@{\,}c@{\,}l}
\dis\bet_k=\sum_{n=k}^\infty e^{-\tet^n}
&=&\dis e^{-\tet^k}\sum_{n=k}^\infty e^{-(\tet^n-\tet^k)}\\[5pt]
&=&\dis e^{-\tet^k}\Big(1+\sum_{m=1}^\infty
e^{-(\tet^{k+m}-\tet^k)}\Big)\\[5pt]
&=&\dis e^{-\tet^k}\Big(1+\sum_{m=1}^\infty
e^{-\tet^k(\tet^m-1)}\Big),
\ea\]
where
\[
1\leq\Big(1+\sum_{m=1}^\infty e^{-\tet^k(\tet^m-1)}\Big)
\leq 1+\sum_{m=1}^\infty e^{-\tet(\tet^m-1)}=:K_\tet.
\]
Thus, $\al_k\leq\bet\leq K_\tet\al_k$ for some finite constant
$K_\tet$ that does not depend on $k$. It follows that
$\log(\al_k)\leq\log(\bet_k)\leq\log(\al_k)+\log(K_\tet)$, hence
\[\ba{l}
\dis\liminf_{k\to\infty}N^{-k}\log(\bet_k)
=\liminf_{k\to\infty}N^{-k}\log(\al_k)\\[5pt]
\dis\quad=\liminf_{k\to\infty}N^{-k}\log(\al_k)
=-\lim_{k\to\infty}\big(\ffrac{\tet}{N}\big)^k=-\infty
\ea\]
for $\tet>N$.}

\subsection{Discussion and outline}\label{S:discus}

\subsection*{Motivation}

Population dynamical models (but not contact processes) on the hierarchical
group have been studied before in e.g.\ \cite{SF83,DG93,Daw00}. The contact
process on the hierarchical group which is the subject of the present paper
may be used to model the spread of an infection in a spatially clustered
population. Taking humans as an example, we may think of a site
$(i_0,i_1,\ldots)$ as an address, where $i_0$ is the house number, $i_1$ the
street, $i_2$ the town, $i_3$ the state and so on. In this example, sites at
hierarchical distance less or equal than $1$, $2$, or $3$ from a given site
are addresses that are in the same street, town, or state, respectively. In
case the $\al_k$ are rapidly decaying, our model describes a situation where
infections between large `blocks' of sites, such as towns or states, are rare,
hence the infection has to overcome certain `bottlenecks' in order to spread
and survive in the long run. In this context, we note that another model
exhibiting such bottlenecks is the one-dimensional contact process in a random
environment with fixed, i.i.d. infection rates between neighboring sites, see
e.g.\ \cite{Lig92}. Another motivation to study contact processes on the
hierarchical group is that they may potentially be used to estimate contact
processes on other lattices from below, including long-range processes on $\Z$
(compare \cite{Dys69}).

Apart from possible applications in population biology, we believe our results
are interesting from a more theoretical point of view because of the way we
prove Theorem~\ref{T:main}~(b). Finding upper bounds on the critical recovery
rate of a contact process is generally easier than finding lower bounds. In
line with this, the proof of Theorem~\ref{T:main}~(a) is rather simple, but
part~(b) is much more involved. In fact, as we explain below, there are only a
few known techniques for proving survival of contact processes in `low'
dimensions, and none seems to work well in our setting. The technique we
finally invented is in its essence a renormalization argument. As such, it is
interesting in the more general program of finding rigorous renormalization
techniques for interacting particle systems.

\subsection*{Renormalization}

It has been recognized long ago that the hierarchical group is especially
suitable for renormalization arguments. There exists an extensive literature
on the Ising model on hierarchical lattices (see, e.g.,
\cite{Dys69,BM87,HHW01}). Moreover, the self-avoiding random walk on a
hierarchical group with `effective' dimension four is treated in \cite{BEI92},
while linearly interacting diffusions and near-critical percolation on the
hierarchical group have been considered in \cite{DG93} and \cite{DG06},
respectively. In these last two papers, in order to get rigorous results, the
authors take a `local mean field limit', meaning that they send the freedom
$N$ of the hierarchical group to infinity and rescale to get nontrivial
limits.

The intuitive idea behind our proof of Theorem~\ref{T:main}~(b) is easily
explained. For given $i=(i_0,i_1,\ldots)\in\om_N$, set
\be
B_i:=\big\{(j,i_0,i_1,\ldots)\in\om_N:j\in\{0,\ldots,N-1\}\big\}
\ee
Then $(B_i)_{i\in\om_N}$ is a collection of blocks $B_i\sub\om_N$,
each $B_i$ containing $N$ sites at distance $1$ from each other. We
would like to consider $B_i$ as a single site in a `renormalized'
lattice, such that $B_i$ can be either infected or healthy. Indeed, if
$N$ is large and $\al_1>\de$, then it can be shown that there exists a
`metastable' state on $B_i$ in which roughly a
$(1-\de/\al_1)$-fraction of the sites is infected, and that
transitions from this metastable state to the all-healthy state are
fast and happen rarely. Thus, as long as $\de/\al_1$ is sufficiently
small, we expect our `renormalized' blocks $B_i$ to behave effectively
as a single site, with an effective `renormalized' recovery rate
$\ti\de$ that is much smaller than the original $\de$. Iterating this
procedure, we expect the system to be more and more stable as we move
up the spatial scale, until, in the limit, we never get extinct.

It may well be that this intuition can be made rigorous in a precise way for a
suitably chosen model in the local mean field limit $N\to\infty$, in the
spirit of \cite{DG93,DG06}. Our motivation, however, was to prove results for
fixed $N$. The problem with renormalization-style arguments for fixed $N$ is
that in this case, one is forced to give exact bounds on how stable the
`metastable' state on $B_i$ is, and how fast transitions between this state
and the all-healthy state are. Moreover, these bounds must be translated into
similar bounds on a renormalized lattice, in a way that can be iterated. If
one tries to do this in a straightforward manner this soon becomes very messy
and technical.

The solution we found for this problem is a technique the second author
learned about from a talk by Tom Kurtz on the look-down construction for
Fleming-Viot processes \cite{DK96,DK99} and that originates from Rogers and
Pitman \cite{RP81}. Basically, this is a technique for adding structure to a
Markov process $X$, such that if in the enriched process $(X,Y)$, one forgets
the added structure $Y$, one obtains back the original process $X$. An
interesting feature of this technique is that in the enriched process $(X,Y)$,
the process $X$ is in general not an autonomous Markov process, i.e., the
dynamics of $X$ depend on $Y$. In practice, we will set up a coupling between
a contact process $X$ on the hierarchical group $\om_2$ with freedom $2$, and
an `added-on' process $\ti Y$ that lives on a renormalized lattice and that is
almost a contact process itself. In particular, $\ti Y$ can be stochastically
estimated from below by a contact process $Y$, which is sufficient for our
purposes. For a more detailed discussion of our methods, we refer the reader
to Sections~\ref{S:coup}--~\ref{S:coupdis} below.

\subsection*{Survival}

Since the direct aim of our renormalization argument is to prove survival, we
conclude this section with a discussion of how survival is proved for contact
processes on other lattices. Since most of the literature deals with
nearest-neighbor processes on graphs, for which there is just a single
infection rate, it has become customary to fix the recovery rate to 1,
consider the infection rate as a variable, and prove upper bounds on the
critical infection rate. By a trivial rescaling of time, we may instead fix
the infection rate and vary the recovery rate, hence any upper bound on the
critical infection rate in the traditional setting can be translated into a
lower bound on the critical recovery rate in our setting.

If $\La$ is an infinite (connected, undirected) graph, then it is always
possible to embed a copy of $\Z$ in $\La$, hence the problem can be reduced to
proving survival of the nearest-neighbor contact process on $\Z$. (It is often
possible to do better than just embedding copy of $\Z$ in $\La$, see
\cite[Thm~VI.4.1]{Lig85}.)

For the nearest-neighbor contact process on $\Z$, we are aware of two
independent proofs that $\de_{\rm c}>0$. If the recovery rate $\de$ is
sufficiently small, then it is not hard to set up a comparison between the
contact process on $\Z$ and oriented percolation on $\Z\times\Z_+$, with a
percolation parameter $p$ close to one. The problem can then be reduced to
showing that $p_{\rm c}<1$ for oriented percolation on $\Z\times\Z_+$, which
is known to follow from a Peierls argument (see \cite[Chapter~5]{Dur88}).

An independent aproach for proving survival of the nearest-neighbor
contact process on $\Z$, which gives a better bound on the critical
value, is the method of Holley and Liggett \cite{HL78} (see also
\cite[Section~IV.1]{Lig85}). Their basic observation is that if there
exists a translation invariant probability law on $\{0,1\}^\Z$ such
that the process $X$ started in this initial law satisfies
\be\label{upmeas}
\dif{t}\P[\exists i\in A\mbox{ s.t.\ }X_t(i)=1]\big|_{t=0}\geq 0
\ee
for all finite $A\sub\La$, then by duality $\P[X_0(0)=1]$ gives a lower bound
on the survival probability of the process started with a single infected
site. Holley and Liggett then explicitly construct a renewal measure that
solves (\ref{upmeas}). Their method has been refined in \cite{Lig95}, leading
to the best rigorous lower bound on $\de_{\rm c}$ available to date.

For lattices different from $\Z$, there exist other, independent methods for
obtaining lower bounds on $\de_{\rm c}$. On $\Z^2$, one may use comparison with
a stochastic Ising model. On $\Z^d$ with $d\geq 3$, one may use comparison
with certain linear systems; this method gives the sharpest known bounds in high
dimensions. (For both these techniques, see \cite[Section~VI.4]{Lig85}.) For
processes on trees, there is a very simple lower bound on $\de_{\rm c}$
resulting from a supermartingale argument (see \cite[Thm~I.4.1]{Lig99}).

In general, one can say that proving survival for contact processes is easier
in higher dimensions. In this context, returning to the hierarchical group, we
mention the following fact. Let $\xi=(\xi_t)_{t\geq 0}$ be a random walk on
$\om_N$ that jumps from a point $i$ to $j$ with rates $a(i,j)$ as in
(\ref{adef}), with
\be\label{ddim}
\al_k=N^{-k(2/d)}\qquad(k\geq 1),
\ee
where $d>0$ is some real constant. Then it can be shown that
\be\label{d2}
P^0[\xi_t=0]\sim t^{-d/2}\,\phi(\log t)\qquad\mbox{as}\quad t\to\infty,
\ee
where $\phi$ is a positive, periodic, continuous real function and $f(t)\sim
g(t)$ means that $f(t)/g(t)\to 1$. Thus, if $d$ is an integer, then such a
random walk is similar to a usual short-range random walk on $\Z^d$. (Indeed,
this is more or less to how Brydges, Evans and Imbrie construct a hierarchical
group with `effective' dimension four in \cite{BEI92}, while the scaling in
\cite{DG93} is chosen so as to mimic the critical dimension for linear
systems, which is two.) Note that in particular, $\xi$ is recurrent if and
only if $d\leq 2$.

\detail{In the note RW.tex it is shown that (\ref{d2}) holds provided that
\[
a(i,j)=\sum_{k=|i-j|}^\infty r_kN^{-k}\quad\mbox{with}\quad
r_k=(1-\tet)\tet^{k-1}
\quad\mbox{and}\quad\tet=N^{-2/d}.
\]
Thus, we need 
\[\ba{r@{\,}c@{\,}l}
\dis a(i,j)&=&\dis\sum_{k=|i-j|}^\infty(1-N^{-2/d})(N^{-2/d})^{k-1}N^{-k}\\[5pt]
&=&\dis(N^{2/d}-1)\sum_{k=|i-j|}^\infty N^{-(1+2/d)k}\\[5pt]
&=&\dis(N^{2/d}-1)N^{-(1+2/d)|i-j|}\sum_{k=0}^\infty N^{-(1+2/d)k}\\[5pt]
&=&\dis\frac{N^{2/d}-1}{1-N^{-(1+2/d)}}N^{-\frac{2}{d}|i-j|}N^{-|i-j|}
\ea\]
Neglecting an overall constant which can be removed by rescaling time, we
arrive at (\ref{ddim}).}

These observations are relevant when we consider comparison with linear
systems as a method to prove survival of contact processes on the hierarchical
group. Indeed, since this technique depends on the transience of the
underlying random walk, for processes with rates $\al_k$ as in (\ref{ddim}),
it seems this technique can only work if $d>2$. Note that our
Theorem~\ref{T:main}~(b) shows that $\de_{\rm c}>0$ for any
$d>0$, and in fact for processes with much faster decaying rates.

If we forget about other `high-dimensional' techniques, this leaves us with
two known techniques for establishing lower bounds on the critical recovery
rate that might be succesful on the hierarchical group: comparison with
oriented percolation plus a Peierls argument, or the method of Holley and
Liggett.

It is not hard to set up a comparison between a contact process on
$\om_N$ and some form of oriented percolation on $\om_N\times\Z_+$
(with a whole set of percolation parameters $p_k$ depending on the
hierarchical distance), but this only moves the problem to proving
that the latter percolates if the $p_k$ are sufficiently large. Since
long-range infections are essential for survival, it is not obvious,
and seems rather difficult, to define suitable contours which could
then be counted and estimated in a Peierls argument.

While we did not spend much time investigating oriented percolation on
$\om_N\times\Z_+$, we did spend a considerable amount of effort trying to
adapt the method of Holley and Liggett. As explained in \cite{Lig95}, the
renewal measure of Holley and Liggett may be interpreted as a certain type of
Gibbs measure with the property that in (\ref{upmeas}), equality holds if $A$
is an interval. The difficult part of the proof is then to show that this
equality for intervals extends to an inequality for general subsets
$A\sub\La$. For the hierarchical group $\om_N$, it is not hard to dream up a
good analogue of Liggett's Gibbs measures and to show that (\ref{upmeas}) may
be satisfied with equality for certain special sets. (In fact, we used blocks
of sites within a given hierarchical distance of each other.) We were not
able, however, to carry out the difficult step in the argument, which is to
extend the equality in (\ref{upmeas}) for special $A$ (the blocks) to an
inequality for general $A\sub\om_N$. It may be that this method can
be carried out succesfully; our failure to do so is no proof that it cannot be
done.

\subsubsection*{Outline}

After proving Theorem~\ref{T:main}~(a) in Section~\ref{S:ext}, we present and
prove our coupling of contact processes on $\om_2$ in Section~\ref{S:coupling}
below. A more detailed discussion of our coupling can be found in
Sections~\ref{S:coup}--~\ref{S:coupdis} while
Sections~\ref{S:couprf}--\ref{S:twolev} contain proofs. The proof of
Theorem~\ref{T:main}~(b) is given in
Section~\ref{S:surv}. Appendix~\ref{A:coord} contains a simple, but rather
tedious argument needed in Section~\ref{S:red}.

\section{Extinction}\label{S:ext}

\subsection{Some general notation}\label{S:notat}

Fix $N\geq 2$ and let $\om=\om_N$ denote the hierarchical group with freedom
$N$. We introduce contact processes whose state spaces are finite analogues of
$\om$. For $n\geq 1$, set
\be
\om^n:=\big\{i=(i_0,\ldots,i_{n-1}):i_k\in\{0,\ldots,N-1\}\big\}
\ee
and
\be
\om^0:=\{(\emptyset)\},
\ee
where $(\emptyset)$ denotes the empty sequence. We equip $\om^n$ with
componentwise addition modulo $N$. For $m,n\geq 0$, we define the
concatenation $i\circ j\in\om^{m+n}$ of elements $i\in\om^m$ and
$j\in\om^n$ by
\be
i\circ j:=(i_0,\ldots,i_{m-1},j_0,\ldots,j_{n-1}).
\ee
Given $0\leq m\leq n$, by definition, the $m$-block in $\om^n$ with
index $j\in\om^{n-m}$ is the set
\be\label{block}
B_m(j):=\{i\circ j:i\in\om^m\}
\qquad(j\in\om^{n-m},\ 0\leq m\leq n).
\ee
We define the set of spin configurations on $\om^n$ by
\be
S_n:=\{0,1\}^{\om^n}=\big\{x=(x(i))_{i\in\om^n}:x(i)\in\{0,1\}\big\}
\qquad(n\geq 0).
\ee
Note that $\om^0$ is a set containing one element and therefore
$S_0=\{0,1\}$. For $0\leq m\leq n$, $i\in\om^{n-m}$, and $x\in S_n$,
we define $x_i\in S_m$ by
\be\label{xij}
x_i(j):=x(j\circ i)\qquad(i\in\om^{n-m},\ j\in\om^m,\ x\in S_n,\ 0\leq m\leq n).
\ee
Note that $x_i$ describes what the spin configuration $x$
looks like on the $m$-block with index~$i$.

For $i\in\om^n$, we define $|i|$ as in (\ref{norm}) with $|i|:=n$ if
$i_{n-1}\neq 0$. For given $\de>0$ and nonnegative constants
$\al_1,\ldots,\al_n$, we define infection rates $a(i,j)$ on $\om^n$ as
in (\ref{adef}), and we call the contact process with
these infection rates and with recovery rate $\de$ the
{\em $(\de,\al_1,\ldots,\al_n)$-contact process}.

\subsection{Extinction}

{\bf Proof of Theorem~\ref{T:main}~(a)} For $n\geq 0$, let $X^{(n)}$ be the
$(\de,\al_1,\ldots,\al_n)$-contact process, and set
\be
l(n):=\E^{\de_0}\big[\inf\{t\geq 0:X^{(n)}_t=\un 0\}\big]\qquad(n\geq 0),
\ee
where $\E^{\de_0}$ denotes expectation with respect to the law of the process
started in $\de_0$ (compare (\ref{Pde}) and note that $0$ now denotes the origin
$0=(0,\ldots,0)\in\om^n$). We will estimate $l(n)$ by a very crude
argument. By a simple rescaling of time, we may assume that the constant $|a|$
in (\ref{summ}) satisfies $|a|=1$. By an obvious coupling, it follows that
$X^{(n)}$ may be stochastically bounded from above by a process $\ti X^{(n)}$
in $S_n=\{0,1\}^{\om^n}$ where sites jump independently of each other from $0$
to $1$ with rate $1$ and from $1$ to $0$ with rate $\de$. Obviously, the
process $\ti X^{(n)}$ has a unique equilibrium law, which is of product form,
and if $\ti X^{(n)}_\infty$ denotes a random variable distributed according to
this law, then
\be\label{prod}
\P[\ti X^{(n)}_\infty=\un 0]=\Big(\frac{\de}{1+\de}\Big)^{N^n}.
\ee
On the other hand, since the Markov process $\ti X^{(n)}$ stays on average a
time $(N^n)^{-1}$ in the state $\un 0$ every time it gets there, one has
\be\label{stay}
\P[\ti X^{(n)}_\infty=\un 0]
=\frac{N^{-n}}{\ti l(n)+N^{-n}}
=\frac{1}{1+N^n\ti l(n)},
\ee
where
\be
\ti l(n):=\E^{\de_0}\big[\inf\{t\geq 0:\ti X^{(n)}_t=\un 0\}\big]\qquad(n\geq 0).
\ee
Solving $\ti l(n)$ from (\ref{prod}) and (\ref{stay}) and comparing with
$l(n)$, we find that
\be\label{lnbd}
l(n)\leq\ti l(n)=N^{-n}\big((1+\de^{-1})^{N^n}-1\big)
\leq N^{-n}(1+\de^{-1})^{N^n}.
\ee

Now consider our original contact process on the (infinite) hierarchical group
$\om_N$. We may stochastically estimate this process from above by a process
where infections over a hierarchical distance $>n$ yield infections of a new
type, in such a way that infections of different types do not interact with
each other (in particular, sites may be infected with infections of
more than one type). Thus, in our new process, each type evolves as a
$(\de,\al_1,\ldots,\al_n)$-contact process in some $n$-block, and in addition,
for each $k>n$, each site that is infected with this type establishes with
rate $\al_kN^{-k}(N^k-N^{k-1})$ another type at a uniformly chosen site in
some uniformly chosen $n$-block at hierarchical distance $k$. Since at any
point in time there are at most $N^n$ infected sites of a given type, and each
type exists for an expected time of length $l(n)$, it follows that the
expected number of new types created by a type during its lifetime is bounded
{f}rom above by
\be
N^nl(n)(1-N^{-1})\sum_{k=n+1}^\infty\al_k.
\ee
In view of (\ref{lnbd}) and the definition of $\bet_n$, we may estimate this
quantity from above by
\be\label{offspring}
(1-N^{-1})(1+\de^{-1})^{N^n}\bet_{n+1}.
\ee
If, for some $n\geq 1$, this quantity is less than 1, then types create new
types according to a subcritical branching process, hence a.s.\ at most
finitely many types are created at all time, hence our
contact process dies out. Taking logarithms and dividing by $N^n$, we see that
for all $\de>0$ there exists an $n\geq 1$ such that the quantity in
(\ref{offspring}) is less than one, provided that
\be
\liminf_{n\to\infty}N^{-n}
\log\big((1-N^{-1})(1+\de^{-1})^{N^n}\bet_{n+1}\big)<0
\qquad\forall\de>0.
\ee
This is equivalent to 
\be
\log(1+\de^{-1})+\liminf_{n\to\infty}N^{-n}\log(\bet_{n+1})<0
\qquad\forall\de>0,
\ee
which is in turn equivalent to (\ref{loglim}).\qed

\section{Coupling of contact processes}\label{S:coupling}

\subsection{A coupling}\label{S:coup}

Throughout this section, we fix $N=2$ and consider finite
$(\de,\al_1,\ldots,\al_n)$-contact processes on $\om^n$ as defined in
Section~\ref{S:notat}. We will prove the following result.
\bp{\bf(Coupling of contact process)}\label{P:coup}
Let $n\geq 1$, $\de>0$, and $\al_1,\ldots,\al_n\geq 0$. Let
$X=(X_t)_{t\geq 0}$ be the $(\de,\al_1,\ldots,\al_n)$-contact process
started in any initial law. Set $\de':=2\xi\de$ and
$\al'_k:=\frac{1}{2}\al_{k+1}$ $(k=1,\ldots,n-1)$, where $\xi=f(\al_1/\de)$
and $f$ denotes the function
\be\label{fdef}
f(r):=\ga-\sqrt{\ga^2-\ffrac{1}{2}}\quad\mbox{with}\quad
\ga:=\ffrac{1}{4}\big(3+\ffrac{1}{2}r\big)\qquad(r\geq 0).
\ee
Then $X$ can be coupled to a process $(\ti Y,Y)$ such that $(X_t,\ti
Y_t)_{t\geq 0}$ is a Markov process, $(Y_t)_{t\geq 0}$ is a
$(\de',\al'_1,\ldots,\al'_{n-1})$-contact process, $\ti Y_0=Y_0$, $\ti
Y_t\geq Y_t$ for all $t\geq 0$, and
\be\label{cond1}
\P\big[\ti Y_t=y\,\big|\,(X_s)_{0\leq s\leq t}\big]=P(X_t,y)\quad{\rm a.s.}
\qquad(t\geq 0,\ y\in S_{n-1}),
\ee
where $P$ is the probability kernel from $S_n$ to $S_{n-1}$ defined by
(recall (\ref{xij}))
\be\label{Pdef1}
P(x,y):=\prod_{i\in\om^{n-1}}p(x_i,y(i))\qquad(x\in S_n,\ y\in S_{n-1}),
\ee
with
\be\label{Pdef2}
\left(\ba{cc}
p(00,0)&p(00,1)\\
p(01,0)&p(01,1)\\
p(10,0)&p(10,1)\\
p(11,0)&p(11,1)\\
\ea\right) 
:=
\left(\ba{cc}
1&0\\
\xi&1-\xi\\
\xi&1-\xi\\
0&1\\
\ea\right).
\ee
\ep

\begin{figure}
\begin{center}
\plaat{4cm}{\includegraphics[width=10cm]{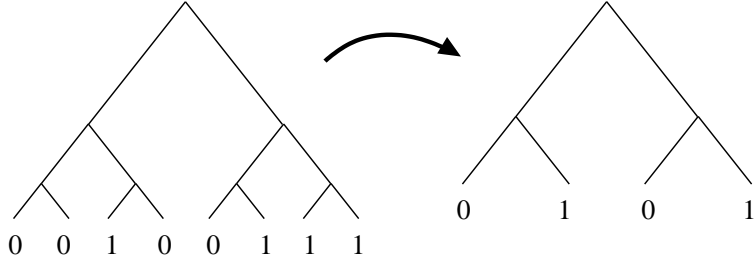}}
\caption[Coupling of $X_t$ and $\ti Y_t$.]
{Coupling of $X_t$ (left) and $\ti Y_t$ (right). The conditional probability
  of the transition depicted here is $\xi(1-\xi)$.}
\label{fig:renorm}
\end{center}
\end{figure}

The coupling in Proposition~\ref{P:coup} achieves the intuitive aim explained
in Section~\ref{S:discus}, namely, to view blocks, consisting of two sites at
distance one from each other, as single sites in a `renormalized' lattice,
which can either be infected or healthy. Indeed, (\ref{cond1}) says that the
conditional law of $\ti Y_t$ given $X_t$ has the following description. First,
we group the sites of $X_t$ into blocks, each consisting of two sites at
distance one from each other. Then, independently for each block, if the
configuration in such a block is $00$ (resp.\ $11$), then we let the
corresponding single site in $\ti Y_t$ be healthy (resp.\ infected), while if
the configuration is $01$ or $10$, then we let the corresponding site in
$\ti Y_t$ be healthy with probability $\xi$ and infected with probability
$1-\xi$. This stochastic rule is demonstrated in Figure~\ref{fig:renorm}. The
transition there has probability $\xi(1-\xi)$ and sites in $\om^3$
and $\om^2$ are depicted as leaves of a binary tree.

It is interesting that a stochastic rule for deciding whether a block is
healthy or infected seems to work better than a deterministic rule. We will
choose the function $p(\,\cdot\,,1)$ in (\ref{Pdef2}) in such a way that this
is the leading eigenfunction of a one-level $(\de,\al_1)$-contact
process. Thus, our methods combine some elements of spectral analysis with
probabilistic coupling tools.

The next lemma (which is proved in Section~\ref{S:addef} below) lists some
elementary properties of the function $f$ defined in (\ref{fdef}).
\bl{\bf(The function $f$)}\label{L:gamma}
The function $f$ defined in (\ref{fdef}) is decreasing on $[0,\infty)$ and
satisfies $f(0)=\frac{1}{2}$ and
\be
f(r)=2r^{-1}+O(r^{-2})\qquad\mbox{as}\quad r\to\infty.
\ee
\el

\subsection{Markov processes with added structure}

In Proposition~\ref{P:coup}, the coupling between the processes $X$ and $\ti
Y$ is of a special kind. There exist general results that tell us how to
construct processes with conditional probabilities as in (\ref{cond1}), such
that in addition $(X_t)_{t\geq 0}$, on its own, is a Markov process. In the
present section, we formulate one such result, which will then be used to
construct the coupling in Proposition~\ref{P:coup}.

Let $S,S'$ be finite sets and set $\hat S:=S\times S'$. Let
$(X,Y)=(X_t,Y_t)_{t\geq 0}$ be a Markov process with state space $\hat
S$ and generator $\hat G$. For each $y\in S'$ (resp.\ $x\in S$), we
define an operator $G_y:\R^S\to\R^S$ (resp.\ $G'_x:\R^{S'}\to\R^{S'}$)
by
\be\ba{rclll}
\dis G_yf(x)&:=&\dis\hat G\ov f(x,y)\quad&\dis\mbox{where}\quad
\ov f(x,y):=f(x)\quad&\dis(x\in S,\ y\in S',\ f\in\R^S),\\[5pt]
\dis G'_xf(y)&:=&\dis\hat G\ov f(x,y)\quad&\dis\mbox{where}\quad
\ov f(x,y):=f(y)\quad&\dis(x\in S,\ y\in S',\ f\in\R^{S'}).
\ec

We say that {\em $X$ evolves according to the generator $G_y$ while $Y=y$}
(resp.\ $Y$ evolves according to the generator $G'_x$ while $X=x$). In
particular, if $G_y$ does not depend on $y$, i.e., if $G_y=G$ $(y\in S')$ for
some operator $G:\R^S\to\R^S$, then we say that $X$ is an {\em autonomous
  Markov process} with generator $G$.  This is equivalent to the statement
that for {\em every} initial law of the joint process $(X,Y)$, the process
$X$, on its own, is the Markov process with generator $G$. The next
proposition (which will be proved in Section~\ref{S:couprf}) demonstrates that
even when $X$ is not autonomous, it may happen that there exists an operator
$G$ such that for certain {\em special} initial laws of the joint process
$(X,Y)$, the process $X$, on its own, is the Markov process with
generator~$G$.  It seems that Rogers and Pitman \cite{RP81} were the first
who noticed this phenomenon. We will prove the proposition below by elaborating
on their result.
\bp{\bf(Markov process with added structure)}\label{P:add}
Let $X$ be a continuous-time Markov process with finite state space
$S$ and generator $G$. Let $S'$ be a finite set, let $P$ be a
probability kernel from $S$ to $S'$, and let $(G'_x)_{x\in S}$ be
a collection of generators of $S'$-valued Markov processes.
Define an operator $\ov G:\R^{S'}\to\R^{S\times S'}$ by
\be\label{ovGdef}
\ov Gf(x,y):=G'_xf(y)\qquad(x\in S,\ y\in S',\ f\in\R^{S'}),
\ee
and define $P:\R^{S'}\to\R^S$ and $\ov P:\R^{S\times S'}\to\R^S$ by
\be\label{Pdef}
Pf(x):=\sum_{y\in S'}P(x,y)f(y)
\quad\mbox{and}\quad
\ov Pf(x):=\sum_{y\in S'}P(x,y)f(x,y).
\ee
Assume that
\be\label{assum}
GPf=\ov P\,\ov G f\qquad(f\in\R^{S'}).
\ee
Then $X$ can be coupled to an $S'$-valued process $Y$ such that
$(X,Y)=(X_t,Y_t)_{t\geq 0}$ is a Markov process with state space
$S\times S'$, the process $Y$ evolves according to the generator
$G'_x$ while $X=x$, and
\be\label{cond}
\P\big[Y_t=y\,\big|\,(X_s)_{0\leq s\leq t}\big]
=P(X_t,y)\quad{\rm a.s.}\qquad(t\geq 0,\ y\in S').
\ee
\ep
\med

\noi
{\bf Remark~1} If $X$ and $Y$ are coupled as in
Proposition~\ref{P:add}, then it is typically not the case that $X$ is
an autonomous Markov process. Nevertheless, the joint Markov process
$(X,Y)$ has the property that if the initial law satisfies
\be\label{condin}
\P\big[Y_0=y\,\big|\,X_0\big]=P(X_0,y)\quad{\rm a.s.}\qquad(t\geq 0,\ y\in S'),
\ee
then $X$, on its own, is the Markov process with generator $G$, and
(\ref{cond}) holds.\med

\noi
{\bf Remark~2} If $X$ and $Y$ are coupled as in
Proposition~\ref{P:add}, then it may happen that $Y$ is an autonomous
Markov process. In this case, we will say that $Y$ is an {\em averaged
Markov process} associated with $X$. In the general case, we will say
that $Y$ is an {\em added-on} process.\med

\subsection{Discussion}\label{S:coupdis}

We mention a few open problems concerning our coupling.

$1^\circ$ Can one modify Proposition~\ref{P:coup} such that $\ti Y=Y$, i.e.,
(in terminology invented in the previous section), for a given
$(\de,\al_1,\ldots,\al_n)$-contact process $X$, can we find a
$(\de',\al'_1,\ldots,\al'_{n-1})$-contact process $Y$ such that $Y$ is an
averaged Markov process of $X$? This would probably involve a kernel $P$ and
constants $\de',\al'_1,\ldots,\al'_{n-1}$ that are more difficult to describe
and less explicit than the ones in Proposition~\ref{P:coup} but would have
great theoretical value, since the resulting map
$(\de,\al_1,\ldots,\al_n)\mapsto(\de',\al'_1,\ldots,\al'_{n-1})$ would
represent a rigorous renormalization transformation.

$2^\circ$ Is it possible to construct a similar coupling as in
Proposition~\ref{P:coup}, but with $\ti Y_t\leq Y_t$? This could potentially
be used to relax condition~(\ref{loglim}).

$3^\circ$ Can one use Proposition~\ref{P:coup} to construct a probability law
on $\{0,1\}^{\om_2}$ that satisfies condition (\ref{upmeas}) of Holley and
Liggett? This would not add much in the line of proving survival (which is
already achieved) but might add to our understanding of the method of Holley
and Liggett, which is rather poor. In particular, in \cite{Lig95} it is shown
that this method may be used to calculate a sequence of approximations of the
critical recovery rate, but beyond the second member of that sequence, there
is no proof that these approximations are lower bounds on $\de_{\rm c}$
(though they are conjectured to be so).

$4^\circ$ Is it possible to make the methods of the present paper work on $\Z$
instead of $\om_2$?  At first sight, it seems that the hierarchical structure
of $\om_2$ is essential to Proposition~\ref{P:coup}. However, when we think of
the latter as `forgetting the fast modes of the spectrum', something may be
possible. Any link between Proposition~\ref{P:coup} and the method of Holley
and Liggett might also provide a clue.

\subsection{Added-on processes}\label{S:couprf}

{\bf Proof of Proposition~\ref{P:add}} We adopt the convention that sums over
$x,x',x''$ always run over $S$ and sums over $y,y',y''$ always run over
$S'$. Write
\bc
\dis Gf(x)&=&\dis\sum_{x'}r(x,x')\big(f(x')-f(x)\big),\\[5pt]
\dis G'_xf(y)&=&\dis\sum_{y'}r'_x(y,y')\big(f(y')-f(y)\big),
\ec
where $r(x,x')$ (resp.\ $r'_x(y,y')$) denotes the rate at which the Markov
process with generator $G$ (resp.\ $G'_x$) jumps from a state $x$ to
a state $x'$ (resp.\ from $y$ to $y'$).

Set $\hat S:=S\times S'$. We let $(X,Y)=(X_t,Y_t)_{t\geq 0}$
be the Markov process in $\hat S$ started in an initial law satisfying
(\ref{condin}), with generator $\hat G$ defined by
\bc\label{hatG2}
\dis\hat Gf(x,y)
&:=&\dis\sum_{x'}t_y(x,x')\big(f(x',y)-f(x,y)\big)\\[5pt]
&&\dis+\sum_{y':\,P(x,y')>0}r'_x(y,y')\big(f(x,y')-f(x,y)\big)\\[5pt]
&&\dis+\sum_{y':\,P(x,y')=0}r'_x(y,y')
\sum_{x'}q_{y'}(x,x')\big(f(x',y')-f(x,y)\big)
\ec
$((x,y)\in\hat S,\ f\in\R^{\hat S})$, where
\be\label{rtq}
t_y(x,x'):=\frac{r(x,x')P(x',y)}{P(x,y)}
\quad\mbox{and}\quad
q_y(x,x'):=\frac{r(x,x')P(x',y)}{\sum_{x''}r(x,x'')P(x'',y)}.
\ee
These formulas are not defined if $P(x,y)=0$ resp.\
$\sum_{x''}r(x,x'')P(x'',y)=0$, so in the first case we define
$t_y(x,x')$, in some arbitrary way, while in the second case we choose
for $q_{y'}(x,\,\cdot\,)$ some arbitrary probability distribution on
$S$. In any case, it will be true that
\be\label{anycase}
\Big(\sum_{x''}r(x,x'')P(x'',y)\Big)q_y(x,x')=r(x,x')P(x',y),
\ee
since the right-hand side of this equation is zero if
$\sum_{x''}r(x,x'')P(x'',y)=0$.

Formula (\ref{hatG2}) says that the process $(X,Y)$ jumps from a state $(x,y)$
to a state $(x',y)$ with rate $t_y(x,x')$. In addition, while $X$ is in the
state $x$, the process $Y$ jumps from the state $y$ to the state $y'$ with
rate $r'_x(y,y')$. During such a jump, if $P(x,y')=0$, then the process $X$
does nothing but if $P(x,y')=0$, then the process $X$ jumps at the same time
to a state $x'$ chosen according to the probability kernel $q_{y'}(x,x')$. In
particular, these rules say that the process $Y$ evolves according to the
generator $G'_x$ while $X=x$.

It is known that (\ref{cond}) holds for the process $(X,Y)$ started
in any initial law satisfying (\ref{condin}), provided that
\be\label{commute}
G\ov Pf=\ov P\hat Gf\qquad(f\in\R^{\hat S}).
\ee
The sufficiency of (\ref{commute}) follows, for example, from
\cite[Corollary~3.5]{Kur98}, which is a rather technical statement
about martingale problems. A much less technical version of this
result can be found in \cite{RP81}.

We may rewrite $\hat G$ in the form
\bc
\dis\hat Gf(x,y)
&:=&\dis\sum_{x'}t_y(x,x')\big(f(x',y)-f(x,y)\big)\\[5pt]
&&\dis+\sum_{y'}r'_x(y,y')\big(f(x,y')-f(x,y)\big)\\[5pt]
&&\dis+\sum_{y':P(x,y')=0}r'_x(y,y')
\sum_{x'}q_{y'}(x,x')\big(f(x',y')-f(x,y')\big)\qquad(x\in S,\ y\in R).
\ec
We calculate, remembering the definition of~$t_y(x,x')$, and letting $G(x,x')$
denote the matrix associated with the operator $G$,
\bc\label{een}
\dis G\ov Pf(x)
&=&\dis\sum_{x'}r(x,x')\Big(\sum_yP(x',y)f(x',y)-\sum_yP(x,y)f(x,y)\Big)\\[5pt]
&=&\dis\sum_y\sum_{x'}r(x,x')\big[P(x',y)\big(f(x',y)-f(x,y)\big)
+\big(P(x',y)-P(x,y)\big)f(x,y)\big]\\[5pt]
&=&\dis\sum_yP(x,y)\sum_{x'}t_y(x,x')\big(f(x',y)-f(x,y)\big)\\[5pt]
&&\dis+\sum_{y:\,P(x,y)=0}\sum_{x'}r(x,x')P(x',y)
\big(f(x',y)-f(x,y)\big)\\[5pt]
&&\dis+\sum_y\Big(\sum_{x'}G(x,x')P(x',y)\Big)f(x,y),
\ec
and
\bc\label{twee}
\dis\ov P\hat Gf(x)
&=&\dis\sum_yP(x,y)\sum_{x'}t_y(x,x')\big(f(x',y)-f(x,y)\big)\\[5pt]
&&\dis+\sum_{yy'}P(x,y)r'_x(y,y')\big(f(x,y')-f(x,y)\big)\\[5pt]
&&\dis+\sum_yP(x,y)\sum_{y':\,P(x,y')=0}r'_x(y,y')
\sum_{x'}q_{y'}(x,x')\big(f(x',y')-f(x,y')\big)\\[5pt]
&=&\dis\sum_yP(x,y)\sum_{x'}t_y(x,x')\big(f(x',y)-f(x,y)\big)\\[5pt]
&&\dis+\sum_{y:\,P(x,y)=0}\sum_{x'}\Big(\sum_{y'}P(x,y')r'_x(y',y)\Big)q_y(x,x')
\big(f(x',y)-f(x,y)\big)\\[5pt]
&&\dis+\sum_{y}\Big(\sum_{y'}P(x,y')G_x(y',y)\Big)f(x,y),
\ec
where to get the second equality we have reordered our terms and
relabelled indices. The first terms on the right-hand sides of
(\ref{een}) and (\ref{twee}) are equal while the third terms agree by
(\ref{assum}). Since by (\ref{assum}), for each $x,y$ such that
$P(x,y)=0$, one has
\be
\sum_{y'}P(x,y')r'_x(y',y)=\sum_{y'}P(x,y')G_x(y',y)=\sum_{x''}G(x,x'')P(x'',y)
=\sum_{x''}r(x,x'')P(x'',y),
\ee
we see by (\ref{anycase}) that also the second terms on the right-hand
sides of (\ref{een}) and (\ref{twee}) agree, hence (\ref{commute})
holds.\qed

\noi
{\bf Remark} For fixed $y\in S'$, set $S_y:=\{x\in S:P(x,y)>0\}$ and
consider the operator $\ti G_y$ defined by (compare (\ref{hatG2})--(\ref{rtq}))
\be
\ti G_yf(x):=\sum_{x'\in S_y}t_y(x,x')\big(f(x')-f(x)\big)
\qquad(x\in S_y,\ f\in\R^{S_y}).
\ee
Then $\ti G_y$ is a `compensated $h$-transform' of the operator $G$,
with the function $h(x):=P(x,y)$. Here, if $G$ is the generator of a
Markov process on $S$ and $h$ is a nonnegative function on $S$, then
\be\label{htrafo}
G^hf:=h^{-1}G(hf)-h^{-1}(Gh)f
\ee
defines a generator of a Markov process on the space
$S_h:=\{x:h(x)>0\}$. This sort of transformation has been called a
{\em compensated $h$-transform} in \cite{FS04}. In particular, if $h$
is harmonic, i.e., $Gh=0$, then $G^h$ is the usual $h$-transform of
$G$.

\subsection{Definition of the added-on process}\label{S:addef}

In this section, we prove Proposition~\ref{P:coup}. Our proof depends
on some calculations that will be done in the next three sections. We
wish to construct an $S^{n-1}$-valued added-on process $\ti Y$ on $X$,
such that $\ti Y$ can be stochastically estimated from below by a
$(\de',\al'_1,\ldots,\al'_{n-1})$-contact process. We introduce the notation
\be
x(i,j):=\big(x(i),x(j)\big)\qquad(x\in S_n,\ i,j\in\om^n).
\ee
With this notation, the generator of $X$ can be written as
\bc\label{Gdef}
\dis Gf(x)
&=&\dis\de\sum_{i\in\om^n}1_{\{x(i)=1\}}\big(f(x-\de_i)-f(x)\big)\\[5pt]
&&\dis+\sum_{k=1}^n\al_k\,2^{-k}\!\!\sum\subb{i,j\in\om^n}{|i-j|=k}\!
1_{\{x(i,j)=(0,1)\}}\big(f(x+\de_i)-f(x)\big).
\ec
For any $x\in S_1=\{0,1\}^2$, we write
\be\label{ovx}
\ov x:=\left\{\ba{ll}
00\quad&\mbox{if }x=(0,0),\\
01\quad&\mbox{if }x=(0,1)\mbox{ or }(1,0),\\
11\quad&\mbox{if }x=(1,1),
\ea\right.
\ee
For each $x\in S_n$, we define a generator $G'_x$ of
an $S_{n-1}$-valued Markov process by (recall (\ref{xij}))
\be\ba{l}\label{Gx2}
\dis G'_xf(y)
=\de'\sum_{i\in\om^{n-1}}1_{\{y(i)=1\}}\big(f(y-\de_i)-f(y)\big)\\[5pt]
\dis\quad+\sum_{k=1}^{n-1}\al_{k+1}\,2^{-k}
\!\!\sum\subb{i,j\in\om^{n-1}}{|i-j|=k}\!
\big[a\big(\ov x_i,\ov x_j\big)1_{\{y(i,j)=(0,1)\}}\\[-20pt]
\dis\phantom{\quad+\sum_{k=1}^{n-1}\al_{k+1}\,2^{-k}
\!\!\sum\subb{i,j\in\om^{n-1}}{|i-j|=k}\!\big[}
+b\big(\ov x_i,\ov x_j\big)1_{\{y(i,j)=(0,0)\}}\big]
\big(f(y+\de_i)-f(y)\big),\\[-20pt]
\ec
where $\de'$ is defined as in Proposition~\ref{P:coup} and $a,b$ are
the functions
\be\label{aa}
\left(\ba{ccc}
a(00,00)&a(00,01)&a(00,11)\\
a(01,00)&a(01,01)&a(01,11)\\
a(11,00)&a(11,01)&a(11,11)\\
\ea\right)
=
\left(\ba{ccc}
\ast&1-\xi&2(1-\xi)\\
\ast&\frac{1}{2}&1\\
\ast&\ast&\ast\\
\ea\right)
\ee
and
\be\label{bb}
\left(\ba{ccc}
b(00,00)&b(00,01)&b(00,11)\\
b(01,00)&b(01,01)&b(01,11)\\
b(11,00)&b(11,01)&b(11,11)\\
\ea\right)
=
\left(\ba{ccc}
0&1-\xi&\ast\\
0&\frac{1}{2}&\ast\\
\ast&\ast&\ast\\
\ea\right).
\ee
Here $\xi$ is defined as in Proposition~\ref{P:coup} and the symbol $\ast$
indicates that the definition of $a$ and $b$ in these points is
irrelevant. Indeed, in the next three sections, we will prove the following
fact.
\bl{\bf(Added-on process)}\label{L:add}
Let $G$ be the generator of the $(\de,\al_1,\ldots,\al_n)$-contact
process on $S_n$, let $P$ be the probability kernel from $S_n$ to
$S_{n-1}$ defined in Proposition~\ref{P:coup}, and let
$(G'_x)_{x\in S_n}$ be the generators defined in (\ref{Gx2}), where
the functions $a$ and $b$ are defined as in
(\ref{aa})--(\ref{bb}). Then, no matter how we define $a$ and $b$ in
points indicated with the symbol $\ast$, one has
\be\label{assum3}
GPf=\ov P\,\ov Gf\qquad(f\in\R^{S_{n-1}}),
\ee
where $\ov G,P$, and $\ov P$ are defined as in (\ref{ovGdef})--(\ref{Pdef}).
\el
Based on Lemma~\ref{L:add}, we can now prove Proposition~\ref{P:coup}.\med

\noi
{\bf Proof of Proposition~\ref{P:coup}} By Proposition~\ref{P:add} and
Lemma~\ref{L:add}, we can couple $X$ to an $S_{n-1}$-valued process
$\ti Y$ such that $(X_t,\ti Y_t)_{t\geq 0}$ is a Markov process, $\ti
Y$ evolves according to the generator $G'_x$ while $X=x$, and
(\ref{cond1}) holds.

By Lemma~\ref{L:gamma}, $0<\xi\leq\frac{1}{2}$. It follows that the functions
$a$ and $b$ in (\ref{aa})--(\ref{bb}) satisfy $a\geq\frac{1}{2}$ and $b\geq
0$. From this and (\ref{Gx2}), it is easy to see that $(X,\ti Y)$ can be
coupled to a $(\de',\al'_1,\ldots,\al'_{n-1})$-contact process $Y$, such that
$\ti Y_0=Y_0$ and $\ti Y_t\geq Y_t$ for all $t\geq 0$.\qed

\noi
For completeness, we give here the:\med

\noi
{\bf Proof of Lemma~\ref{L:gamma}} Set
$\xi(\ga):=\ga-\sqrt{\ga^2-\ffrac{1}{2}}$. Then it is straightforward to check
that $\xi(\frac{3}{4})=\frac{1}{2}$. Moreover,
$\dif{\ga}\xi(\ga)=1-\ga(\ga^2-\frac{1}{2})^{-1/2}
=1-(1-\frac{1}{2}\ga^{-2})^{-1/2}<0$ on $[\frac{3}{4},\infty)$, so
$\ga\mapsto\xi(\ga)$ is decreasing on $[\frac{3}{4},\infty)$. Set
$\eps:=\ga^{-1}$. Then $\xi(\eps^{-1})
=\eps^{-1}\big(1-\sqrt{1-\frac{1}{2}\eps^2}\big)$. We observe that
\be\ba{l}
\dis\big(1-\sqrt{1-\ffrac{1}{2}\eps^2}\big)\big|_{\eps=0}=0,\\[5pt]
\dis\dif{\eps}\big(1-\sqrt{1-\ffrac{1}{2}\eps^2}\big)\big|_{\eps=0}
=\ffrac{1}{2}\eps(1-\ffrac{1}{2}\eps^2)^{-1/2}\big|_{\eps=0}=0,\\[5pt]
\dis\diff{\eps}\big(1-\sqrt{1-\ffrac{1}{2}\eps^2}\big)\big|_{\eps=0}
=\dif{\eps}\ffrac{1}{2}\eps(1-\ffrac{1}{2}\eps^2)^{-1/2}\big|_{\eps=0}\\[5pt]
\dis\quad=\dif{\eps}\Big(\ffrac{1}{2}(1-\ffrac{1}{2}\eps^2)^{-1/2}
+\ffrac{1}{4}\eps^2(1-\ffrac{1}{2}\eps^2)^{-3/2}\Big)\big|_{\eps=0}
=\ffrac{1}{2},
\ec
hence $\xi(\eps^{-1})=\eps^{-1}\big(\ffrac{1}{4}\eps^2+O(\eps^3)\big)
=\ffrac{1}{4}\eps+O(\eps^2)$, i.e.,
\be
\xi(\ga)=\ffrac{1}{4}\ga^{-1}+O(\ga^{-2})
\qquad\mbox{as}\quad\ga\to\infty.
\ee
To translate this to the statements in Lemma~\ref{L:gamma}, it suffices to
note that the function $r\mapsto\ga(r):=\frac{1}{4}(2+\frac{1}{2}r)$ is
increasing on $\half$, satisfies $\ga(0)=\ffrac{3}{4}$, and
$\ga(r)=\frac{1}{8}r+O(1)$ as $r\to\infty$.\qed

\subsection{Reduction to a one- and two-level system}\label{S:red}

In this section, we prove Lemma~\ref{L:add}. Our proof is based on two
lemmas which will be proved in the next two sections.\med

\noi
{\bf Proof of Lemma~\ref{L:add}} We start by rewriting the generator
in (\ref{Gdef}) as follows:
\bc\label{Gsplit}
\dis Gf(x)
&=&\dis\de\sum_{i\in\om^{n-1}}\sum_{i'\in\om^1}1_{\{x(i'\circ i)=1\}}
\big(f(x-\de_{i'\circ i})-f(x)\big)\\[15pt]
&&\dis+\al_12^{-1}\sum_{i\in\om^{n-1}}\sum\subb{i',i''\in\om^1}{|i'-i''|=1}
1_{\{x(i'\circ i)=0,\ x(i''\circ i)=1\}}\big(f(x+\de_{i'\circ i})-f(x)\big)
\\[-5pt]
&&\dis+\sum_{k=2}^n\al_k\,2^{-k}\!\!\sum\subb{i,j\in\om^{n-1}}{|i-j|=k-1}\!
\sum_{i',j'\in\om^1}1_{\{x(i'\circ i)=0,\ x(j'\circ j)=1\}}
\big(f(x+\de_{i'\circ i})-f(x)\big)\\
&=&\dis\sum_{i\in\om^{n-1}}R_if(x)
+\sum_{k=1}^{n-1}\al_{k+1}\,2^{-k}
\!\!\sum\subb{i,j\in\om^{n-1}}{|i-j|=k}\!I_{ij}f(x),
\ec
where
\bc
\dis R_if(x)&:=&\dis\de\sum_{i'\in\{0,1\}}1_{\{x(i'\circ i)=1\}}
\big(f(x-\de_{i'\circ i})-f(x)\big)\\
&&\dis+\al_12^{-1}\!\!\sum\subb{i',i''\in\{0,1\}}{i'\neq i''}\!
1_{\{x(i'\circ i)=0,\ x(i''\circ i)=1\}}\big(f(x+\de_{i'\circ i})-f(x)\big),\\
\dis I_{ij}f(x)&:=&\dis 2^{-1}\sum_{i',j'\in\{0,1\}}
1_{\{x(i'\circ i)=0,\ x(j'\circ j)=1\}}\big(f(x+\de_{i'\circ i})-f(x)\big).
\ec
Likewise, we may write the operator in (\ref{Gx2}) as
\be\ba{l}\label{Gxsplit}
\dis G'_xf(y)=\sum_{i\in\om^{n-1}}R'_if(y)+\sum_{k=1}^{n-1}\al_{k+1}\,2^{-k}
\!\!\sum\subb{i,j\in\om^{n-1}}{|i-j|=k}{I'}^x_{ij}f(y),\\[-10pt]
\ec
where
\bc
\dis R'_if(y)&:=&\dis\de'1_{\{y(i)=1\}}\big(f(y-\de_i)-f(y)\big),\\[5pt]
\dis {I'}^x_{ij}f(y)&:=&\dis\big[a\big(\ov x_i,\ov x_j\big)1_{\{y(i)=0,\ y(j)=1\}}
+b\big(\ov x_i,\ov x_j\big)1_{\{y(i)=0,\ y(j)=0\}}\big]
\big(f(y+\de_i)-f(y)\big).
\ec
In view of (\ref{Gsplit}) and (\ref{Gxsplit}), in order to prove
(\ref{assum3}), it suffices to show that
\be\left.\ba{rr@{\,}c@{\,}l}\label{RIi}
{\rm(i)}&\dis R_iPf&=&\dis PR'_if\\[5pt]
{\rm(ii)}&\dis I_{ij}Pf&=&\dis\ov P\,\ov I_{ij}f
\ea\right\}\quad(f\in\R^{S_{n-1}},\ i,j\in\om^{n-1},\ i\neq j),
\ee
where $\ov I_{ij}:\R^{S_{n-1}}\to\R^{S_n\times S_{n-1}}$ is defined
as $\ov I_{ij}f(x,y):={I'}^x_{ij}f(y)$. Note that since $R'_i$ does
not depend on $x$, there is no need to define $\ov R_i$.

The operators $R_i,R'_i,I_{ij}$, and $\ov I_{ij}$ act only on certain
coordinates. In view of this, our problem reduces to a
lower-dimensional one, and (\ref{RIi}) follows from
Lemmas~\ref{L:onelev} and \ref{L:twolev} stated below. It is not
difficult, but notationally cumbersome, to give a formal derivation of
(\ref{RIi}) from Lemmas~\ref{L:onelev} and \ref{L:twolev}. For
completeness, we give this derivation in Appendix~\ref{A:coord}.\qed

\noi
Recall that $S_1=\{0,1\}^{\om^1}$ and
$S_0=\{0,1\}^{\om^0}=\{0,1\}$. Let $\de>0$, $\al_1\geq 0$, and let
$\de'$ be defined as in Proposition~\ref{P:coup}. Let $P$ be the
probability kernel from $S_1$ to $S_0$ defined in
(\ref{Pdef1})--(\ref{Pdef2}) and let $P:\R^{S_0}\to\R^{S_1}$ be
defined as in (\ref{Pdef}). Let $R$ be the generator of the
$(\de,\al_1)$-contact process $X$ on $S_1$ and let $R'$ be the
generator of the $\de'$-contact process $Y$ on $S_0$. The latter is
just the Markov process with state space $\{0,1\}$ that jumps from $1$
to $0$ with rate $\de'$. The next lemma implies that $Y$ is an
averaged Markov process associated with $X$, i.e., $X$ and $Y$ can be
coupled such that (\ref{cond}) holds.
\bl{\bf(One-level system)}\label{L:onelev}
One has
\be\label{onelev}
RPf=PR'f\qquad(f\in\R^{S_0}).
\ee
\el
Formula (\ref{onelev}) implies (\ref{RIi})~(i). We next formulate a
lemma that implies (\ref{RIi})~(ii). 

Let $\de>0$, $\al_1\geq 0$, let $P$ be the probability kernel from
$S_2$ to $S_1$ defined in (\ref{Pdef1})--(\ref{Pdef2}), and let $a,b$
be the functions in (\ref{aa})--(\ref{bb}). We define a generator $I$
of a Markov process in $S_2$ and generators $(I'_x)_{x\in S_2}$ of
Markov processes in $S_1$ by
\bc
\dis If(x)&:=&\dis\ffrac{1}{2}\!\!\sum_{i,j\in\{0,1\}}\!
1_{\{x(i,0)=0,\ x(j,1)=1\}}\big(f(x+\de_{(i,0)})-f(x)\big),\\[5pt]
\dis I'_xf(y)&:=&\dis\big[a(\ov x_0,\ov x_1)1_{\{y(0,1)=(0,1)\}}
+b(\ov x_0,\ov x_1)1_{\{y(0,1)=(0,0)\}}\big]\big(f(y+\de_0)-f(y)\big).
\ec
Our next lemma says that the $(I'_x)_{x\in S_2}$ define an added-on 
Markov process associated with the process with generator $I$.
\bl{\bf(Two-level system)}\label{L:twolev}
One has
\be\label{twolev}
IPf=\ov P\,\ov If\qquad(f\in\R^{S_1}),
\ee
where $P,\ov P$ are defined as in (\ref{Pdef}) and $\ov If(x,y):=I'_xf(y)$
$(x\in S_2,\ y\in S_1)$.
\el

\subsection{The one-level system}

{\bf Proof of Lemma~\ref{L:onelev}} We may write (\ref{onelev}) in the
matrix form
\be
\sum_{x'\in S_1}\sum_{y\in S_0}R(x,x')P(x',y)f(y)
=\sum_{y'\in S_0}\sum_{y\in S_0}P(x,y')R'(y',y)f(y)
\qquad(x\in S_1,\ f\in\R^{S_0}),
\ee
which is equivalent to
\be\label{Rassum}
\sum_{x'\in S_1}R(x,x')p(x',y)
=\sum_{y'\in S_0}p(x,y')R'(y',y)
\qquad(x\in S_1,\ y\in S_0),
\ee
where $p$ is the function in (\ref{Pdef2}). Here
\be
\sum_{x'\in S_1}R(x,x')p(x',y)=Rp(\,\cdot\,,y)(x)
\qquad(x\in\{0,1\}^2,\ y\in\{0,1\}).
\ee
Thus, (\ref{Rassum}) says that $Rp(\,\cdot\,,0)$ and $Rp(\,\cdot\,,1)$
can be written as a linear combination of the functions
$p(\,\cdot\,,0)$ and $p(\,\cdot\,,1)$. It follows that $\Fi:={\rm
span}\{p(\,\cdot\,,0),p(\,\cdot\,,1)\}$ is an invariant subspace of
the operator $R$. Since $p(\,\cdot\,,0)+p(\,\cdot\,,1)=1$, the space
$\Fi$ contains the constant function $1$. We will show that $\Fi$ is
in fact the span of $1$ and one nontrivial eigenfunction of $R$.

We start by noting that by symmetry, the space
$\Hi:=\{f\in\R^{S_1}:f(0,1)=f(1,0)\}$ is invariant under $R$. Since
$(0,0)$ is a trap of the $(\de,\al_1)$-contact process, the space
$\Hi_0:=\{f\in\Hi:f(0,0)=0\}$ is also invariant under $R$; in fact,
$\Hi$ is the span of $\Hi_0$ and the trivial eigenfunction $1$.
In view of this, we look for eigenfunctions of $R$ in $\Hi_0$.
We observe that for $f\in\Hi_0$,
\bc
\dis\left(\ba{c}Rf(0,1)\\ Rf(1,1)\ea\right)
&=&\dis\left(\ba{c}\de\big(0-f(0,1)\big)
+\ffrac{1}{2}\al_1\big(f(1,1)-f(0,1)\big)\\
2\de\big(f(0,1)-f(1,1)\big)\ea\right)\\[15pt]
&=&\dis\left(\ba{cc}
-(\de+\ffrac{1}{2}\al_1)&\ffrac{1}{2}\al_1\\
2\de&-2\de
\ea\right)\left(\ba{c}f(0,1)\\ f(1,1)\ea\right).
\ec
To find the eigenvalues,
we must solve
\be
{\rm det}\left(\ba{cc}
-(\de+\ffrac{1}{2}\al_1)-\la&\ffrac{1}{2}\al_1\\[5pt]
2\de&-2\de-\la
\ea\right)=0,
\ee
which gives
\be\ba{l}
\dis(\de+\ffrac{1}{2}\al_1+\la)(2\de+\la)=\de\al_1\\[3pt]
\dis\desd\quad\la^2+(3\de+\ffrac{1}{2}\al_1)\la+2\de^2=0\\
\dis\desd\quad\Big(\la+\frac{3\de+\ffrac{1}{2}\al_1}{2}\Big)^2
=\Big(\frac{3\de+\ffrac{1}{2}\al_1}{2}\Big)^2-2\de^2\\
\dis\desd\quad\la=-\Big(\frac{3\de+\ffrac{1}{2}\al_1}{2}\Big)
\pm\sqrt{\Big(\frac{3\de+\ffrac{1}{2}\al_1}{2}\Big)^2-2\de^2}\\
\dis\desd\quad\la=-2\de\big(\ga\pm\sqrt{\ga^2-\ffrac{1}{2}}\big),
\ec
where $\ga:=\frac{1}{4}(3+\frac{1}{2}\frac{\al_1}{\de})$ (compare
(\ref{fdef})). In particular, the leading eigenvalue is $\la=-2\de\xi=-\de'$,
where $\xi$ and $\de'$ are defined as in Proposition~\ref{P:coup}. To find the
corresponding eigenfunction, we need to solve
\be\ba{l}
\dis 2\de\big(f(0,1)-f(1,1)\big)=\la f(1,1)\\
\dis\desd\quad 2\de f(0,1)=(2\de+\la)f(1,1)=2\de(1-\xi)f(1,1)\\
\dis\desd\quad f(0,1)=(1-\xi)f(1,1),
\ec
which yields the eigenfunction
\be
\left(\ba{c}f(0,0)\\ f(0,1)\\ f(1,0)\\ f(1,1)\ea\right)
=\left(\ba{c}0\\ 1-\xi \\ 1-\xi \\ 1\ea\right)
=p(\,\cdot\,,1).
\ee
Our calculations so far show that $\Fi:={\rm
span}\{1,p(\,\cdot\,,1)\}={\rm span}\{p(\,\cdot\,,0),p(\,\cdot\,,1)\}$
is an invariant subspace of the operator $R$. It follows that there
exist constants $(R'(y',y))_{y,y'\in\{0,1\}}$ such that
\be
Rp(\,\cdot\,,y)=\sum_{y'\in\{0,1\}}p(\,\cdot\,,y')R'(y',y).
\ee
In fact
\be
\big(Rp(\,\cdot\,,0)\ Rp(\,\cdot\,,1)\big)
=\big(\de'p(\,\cdot\,,1)\ \;-\de'p(\,\cdot\,,1)\big)
=\ba{c}\big(p(\,\cdot\,,0)\ p(\,\cdot\,,1)\big)\\
\phantom{p(\,\cdot\,,0)\ p(\,\cdot\,,1)}\ea
\left(\ba{cc}0&0\\ \de'&-\de'\ea\right),
\ee
hence
\be
\left(\ba{cc}R'(0,0)&R'(0,1)\\ R'(1,0)&R'(1,1)\ea\right)
=\left(\ba{cc}0&0\\ \de'&-\de'\ea\right),
\ee
which we recognize as the generator of a Markov process on $\{0,1\}$
that jumps from $1$ to $0$ with rate $\de'$.\qed

\subsection{The two-level system}\label{S:twolev}

{\bf Proof of Lemma~\ref{L:twolev}} We may write (\ref{twolev}) in the
matrix form
\be
\sum_{x'\in S_2}\sum_{y\in S_1}I(x,x')P(x',y)f(y)
=\sum_{y'\in S_1}\sum_{y\in S_1}P(x,y')I'_x(y',y)f(y)
\qquad(x\in S_2,\ f\in\R^{S_1}),
\ee
which is equivalent to
\be\label{Iassum}
\sum_{x'\in S_2}I(x,x')P(x',y)
=\sum_{y'\in S_1}P(x,y')I'_x(y',y)
\qquad(x\in S_2,\ y\in S_1).
\ee
Here
\be
\sum_{x'\in S_2}I(x,x')P(x',y)=IP(\,\cdot\,,y)(x)\qquad(x\in S_2,\ y\in S_1).
\ee
Note that $S_n=\{0,1\}^{\om^n}=\{0,1\}^{\{0,1\}^n}$ has $2^{2^n}$
elements, so $|S_1|=2^2=4$ and $|S_2|=2^4=16$, hence
$(IP(\,\cdot\,,y)(x))_{x\in S_1,\ y\in S_2}$ is a matrix with
$4\cdot16=64$ entries. Luckily, using symmetry, we can reduce the size
of our problem quite a bit. We start by calculating
\be
P(x,y)=P_y(x_0,x_1)=p(x_0,y(0))p(x_1,y(1))
\ee
for $x_0,x_1,y\in\{(0,0),(0,1),(1,1)\}$. For brevity, we write
$00=(0,0)$, $01=(0,0)$, and $11=(1,1)$. We have
\be\label{P00}
\left(\ba{ccc}
P_{00}(00,00)&P_{00}(00,01)&P_{00}(00,11)\\
P_{00}(01,00)&P_{00}(01,01)&P_{00}(01,11)\\
P_{00}(11,00)&P_{00}(11,01)&P_{00}(11,11)
\ea\right)
=\left(\ba{ccc}
1&\xi&0\\
\xi&\xi^2&0\\
0&0&0
\ea\right),
\ee
\be\label{P01}
\left(\ba{ccc}
P_{01}(00,00)&P_{01}(00,01)&P_{01}(00,11)\\
P_{01}(01,00)&P_{01}(01,01)&P_{01}(01,11)\\
P_{01}(11,00)&P_{01}(11,01)&P_{01}(11,11)
\ea\right)
=\left(\ba{ccc}
0&1-\xi&1\\
0&\xi(1-\xi)&\xi\\
0&0&0
\ea\right),
\ee
and
\be\label{P11}
\left(\ba{ccc}
P_{11}(00,00)&P_{11}(00,01)&P_{11}(00,11)\\
P_{11}(01,00)&P_{11}(01,01)&P_{11}(01,11)\\
P_{11}(11,00)&P_{11}(11,01)&P_{11}(11,11)
\ea\right)
=\left(\ba{ccc}
0&0&0\\
0&(1-\xi)^2&1-\xi\\
0&1-\xi&1
\ea\right).
\ee
Recall the definition of $\ov x$ from (\ref{ovx}). If $(X(t))_{t\geq
0}=(X_0(t),X_1(t))_{t\geq 0}$ is a Markov process in $S_2=S_1\times
S_1$ with generator $I$, then $(\ov X_0(t),\ov X_1(t))_{t\geq 0}$
is a Markov process that jumps with the following rates:
\be\ba{c@{\qquad}c@{\qquad}c}
(00,00)&(00,01)&(00,11)\\
&\Big\downarrow 1&\Big\downarrow 2\\
(01,00)&(01,01)&(01,11)\\
&\Big\downarrow\frac{1}{2}&\Big\downarrow 1\\
(11,00)&(11,01)&(11,11).
\ea\ee
{F}rom this, we see that the functions $IP(\,\cdot\,,y)(x)=IP_y(x)$ are given by
\be\label{IP00}
\left(\ba{ccc}
IP_{00}(00,00)&IP_{00}(00,01)&IP_{00}(00,11)\\
IP_{00}(01,00)&IP_{00}(01,01)&IP_{00}(01,11)\\
IP_{00}(11,00)&IP_{00}(11,01)&IP_{00}(11,11)
\ea\right)
=\left(\ba{ccc}
0&-\xi(1-\xi)&0\\
0&-\frac{1}{2}\xi^2&0\\
0&0&0
\ea\right),
\ee
\be\label{IP01}
\left(\ba{ccc}
IP_{01}(00,00)&IP_{01}(00,01)&IP_{01}(00,11)\\
IP_{01}(01,00)&IP_{01}(01,01)&IP_{01}(01,11)\\
IP_{01}(11,00)&IP_{01}(11,01)&IP_{01}(11,11)
\ea\right)
=\left(\ba{ccc}
0&-(1-\xi)^2&-2(1-\xi)\\
0&-\frac{1}{2}\xi(1-\xi)&-\xi\\
0&0&0
\ea\right),
\ee
and
\be\label{IP11}
\left(\ba{ccc}
IP_{11}(00,00)&IP_{11}(00,01)&IP_{11}(00,11)\\
IP_{11}(01,00)&IP_{11}(01,01)&IP_{11}(01,11)\\
IP_{11}(11,00)&IP_{11}(11,01)&IP_{11}(11,11)
\ea\right)
=\left(\ba{ccc}
0&(1-\xi)^2&2(1-\xi)\\
0&\frac{1}{2}\xi(1-\xi)&\xi\\
0&0&0
\ea\right).
\ee
We wish to express the functions $(IP(\,\cdot\,,y))_{y\in S_1}$ in the
functions $(P(\,\cdot\,,y))_{y\in S_1}$. Unlike in the previous
section, the span of the functions $(P(\,\cdot\,,y))_{y\in S_1}$ is
not invariant under the operator $I$, so we cannot express the
functions $(IP(\,\cdot\,,y))_{y\in S_1}$ as a linear combination of
the functions $(P(\,\cdot\,,y))_{y\in S_1}$. However, we can find
expressions of the form (compare (\ref{Iassum}))
\be
IP(\,\cdot\,,y)(x)
=\sum_{y'\in S_1}P(x,y')I'_x(y',y)
\qquad(x\in S_2,\ y\in S_1),
\ee
where the coefficients $I'_x(y',y)$ do not depend too strongly on
$x$. Solutions to this problem are not unique. The claim of
Lemma~\ref{L:twolev} is that we can choose
\be\ba{l}
\dis\left(\ba{cccc}
I'_x(00,00)&I'_x(00,01)&I'_x(00,10)&I'_x(00,11)\\
I'_x(01,00)&I'_x(01,01)&I'_x(01,10)&I'_x(01,11)\\
I'_x(10,00)&I'_x(10,01)&I'_x(10,10)&I'_x(10,11)\\
I'_x(11,00)&I'_x(11,01)&I'_x(11,10)&I'_x(11,11)
\ea\right)\\[30pt]
\dis\qquad\qquad=\left(\ba{cccc}
-b(\ov x_0,\ov x_1)&0&b(\ov x_0,\ov x_1)&0\\
0&-a(\ov x_0,\ov x_1)&0&a(\ov x_0,\ov x_1)\\
0&0&0&0\\
0&0&0&0
\ea\right),
\ec
where $a,b$ are the functions in (\ref{aa})--(\ref{bb}). Thus, we need
to check that
\be\ba{rr@{\,}c@{\,}l}
{\rm(i)}&\dis IP(\,\cdot\,,00)(x)&=&-b(\ov x_0,\ov x_1)P(x,00),\\
{\rm(ii)}&\dis IP(\,\cdot\,,01)(x)&=&-a(\ov x_0,\ov x_1)P(x,01),\\
{\rm(iii)}&\dis IP(\,\cdot\,,10)(x)&=&b(\ov x_0,\ov x_1)P(x,00),\\
{\rm(iv)}&\dis IP(\,\cdot\,,11)(x)&=&a(\ov x_0,\ov x_1)P(x,01).
\ec
Since $\sum_{y\in S_1}IP(\,\cdot\,,y)=I1=0$, it suffices to check only three
of these equations, say (i), (ii), and (iv). We observe from
(\ref{IP01})--(\ref{IP11}) that $IP(\,\cdot\,,01)=-IP(\,\cdot\,,11)$. In view
of this, it suffices to check only (i) and (ii). By (\ref{aa})--(\ref{bb}),
(\ref{P00})--(\ref{P01}), and (\ref{IP00})--(\ref{IP01}), we need to check
that
\be\label{check00}
\left(\ba{ccc}
0&-\xi(1-\xi)&0\\
0&-\frac{1}{2}\xi^2&0\\
0&0&0
\ea\right)
=-\left(\ba{ccc}
0&1-\xi&\ast\\
0&\frac{1}{2}&\ast\\
\ast&\ast&\ast\\
\ea\right)
\bullet\left(\ba{ccc}
1&\xi&0\\
\xi&\xi^2&0\\
0&0&0
\ea\right)
\ee
and
\be\label{check01}
\left(\ba{ccc}
0&-(1-\xi)^2&-2(1-\xi)\\
0&-\frac{1}{2}\xi(1-\xi)&-\xi\\
0&0&0
\ea\right)
=-\left(\ba{ccc}
\ast&1-\xi&2(1-\xi)\\
\ast&\frac{1}{2}&1\\
\ast&\ast&\ast\\
\ea\right)
\bullet\left(\ba{ccc}
0&1-\xi&1\\
0&\xi(1-\xi)&\xi\\
0&0&0
\ea\right),
\ee
where $\bullet$ denotes the componentwise product of functions and the
symbol $\ast$ indicates that the value of the functions $a$ and $b$ in
these points is irrelevant. We see by inspection that (\ref{check00})
and (\ref{check01}) are satisfied.\qed

\section{Survival}\label{S:surv}

\subsection{Survival bounds}

Until further notice, we continue to study the contact process on the
hierarchical group $\om_N$ with $N=2$ and its finite analogues defined
in Section~\ref{S:coup}. Our proof of Theorem~\ref{T:main}~(b) is based on
the following basic estimate.
\bp{\bf(Survival bound for finite systems)}\label{P:finsurv}
Let $\de>0$ and let $(\al_k)_{k\geq 1}$ be nonnegative constants.
Let $X^{(n)}$ be the $(\de,\al_1,\ldots,\al_n)$-contact process
started in $X^{(n)}_0=\de_0$. Then
\be\label{finbd}
\P^{\de_0}[X^{(n)}_t\neq\un 0]\geq\Big(\prod_{k=0}^{n-1}(1-\xi(k))\Big)
\ex{-\de(n)t}\qquad(t\geq 0),
\ee
where $\de(0):=\de$, $\al_k(0):=\al_k$ $(k\geq 1)$, and we
define inductively, for $n\geq 0$,
\bc\label{induc2}
\dis\de(n+1)&:=&\dis 2\xi(n)\de(n),\\[5pt]
\dis\al_k(n+1)&:=&\dis\ffrac{1}{2}\al_{k+1}(n)\qquad(k\geq 1),
\ec
where $\xi(n):=f(\al_1(n)/\de(n))$ with $f$ as in (\ref{fdef}).
\ep
{\bf Proof} By Lemma~\ref{L:gamma}, one has
$0<\xi(n)\leq\frac{1}{2}$ for all $n\geq 0$. For
$0<\xi\leq\frac{1}{2}$ and $k\geq 1$, let $P_{k,\xi}$ denote the
probability kernel from $S_k$ to $S_{k-1}$ defined in
(\ref{Pdef1})--(\ref{Pdef2}). Let $X^{(n)}$ be the
$(\de,\al_1,\ldots,\al_n)$-contact process. Applying
Proposition~\ref{P:coup} inductively, we can couple $X^{(n)}$ to
processes
\[
\ti X^{(n-1)},X^{(n-1)},\ldots,\ti X^{(0)},X^{(0)}
\]
such that $\ti X^{(n-m)},X^{(n-m)}$ take values in $S_{n-m}$, one has
$\ti X^{(n-m)}_0=X^{(n-m)}_0$, $\ti X^{(n-m)}_t\geq X^{(n-m)}_t$ for
all $t\geq 0$,
\be
\P\big[\ti X^{(n-m-1)}_t=y\,\big|\,X^{(n-m)}_t=x\big]=P_{n-m,\xi(m)}(x,y),
\ee
and the process $X^{(n-m)}$ is a
$(\de(m),\al_1(m),\ldots,\al_n(m))$-contact process.  A little
thinking convinces us that this coupling can be done in a Markovian
way, i.e., in such a way that $(\ti X^{(n-m-1)},X^{(n-m-1)})$ is
conditionally independent of
\[
X^{(n)},(\ti X^{(n-1)},X^{(n-1)}),\ldots,(\ti X^{(n-m+1)},X^{(n-m+1)})
\]
given $(\ti X^{(n-m)},X^{(n-m)})$. By this Markovian property and the
definition of $P_{k,\xi}(x,y)$, if we start $X^{(n)}$ in the initial
state $X^{(n)}_0=\de_0$, then
\be
\P[X^{(n-m)}_0=\de_0]=\prod_{k=0}^{m-1}(1-\xi(k)),
\ee
and $X^{(n-m)}_0=\un 0$ with the remaining probability. Since
$\P[X^{(n-m-1)}_t=\un 0\,|\,X^{(n-m)}_t=\un 0]=1$ for each $m$,
we have
\be
\P[X^{(n)}_t\neq\un 0]\geq\P[X^{(n-m)}_t\neq\un 0]\qquad(0\leq m\leq n).
\ee
In particular, since $X^{(0)}$ is a Markov process in $S_0=\{0,1\}$ that
jumps from $1$ to $0$ with rate $\de(n)$, we observe that
\be
\P[X^{(n)}_t\neq\un 0]\geq\P[X^{(0)}_t\neq\un 0]=e^{-\de(n)t}\,\P[X^{(0)}_0=\de_0]
=e^{-\de(n)t}\prod_{k=0}^{n-1}(1-\xi(k))\qquad(t\geq 0),
\ee
which proves (\ref{finbd}).\qed

\noi
As an immediate corollary to Proposition~\ref{P:finsurv}, we obtain:
\bp{\bf(Survival bound for infinite systems)}\label{P:infsurv}
Let $\de>0$ and let $(\al_k)_{k\geq 1}$ be nonnegative constants
satisfying $\sum_{k=1}^\infty\al_k<\infty$. Let $(\xi(k))_{k\geq 0}$
be defined as in Proposition~\ref{P:finsurv}. Let $X$ be the contact process on
$\om_2$ with infection rates as in (\ref{adef}) and recovery rate
$\de$. Then the process started in $X_0=\de_0$ satisfies
\be\label{infbd}
\P^{\de_0}[X_t\neq\un 0\ \forall t\geq 0]\geq\prod_{k=0}^\infty(1-\xi(k)).
\ee
\ep
{\bf Proof} It is easy to see that the process $X^{(n)}$ in
Proposition~\ref{P:finsurv} and $X$ can be coupled such that
$X^{(n)}_t\leq X_t$ for all $t\geq 0$. Therefore (\ref{infbd}) follows
{f}rom (\ref{finbd}), provided we show that $\de(k)\to 0$ as $k\to\infty$.
In fact, it suffices to prove this under the assumption that
$\prod_{k=0}^\infty(1-\xi(k))>0$, for otherwise (\ref{infbd}) is
trivial. Indeed, $\prod_{k=0}^\infty(1-\xi(k))>0$ implies that $\xi(k)\to 0$
as $k\to\infty$, which by the fact that
\be
\de(n)=\de\prod_{k=0}^{n-1}\big(2\xi(k)\big)
\ee
implies that $\de(k)\to 0$ as $k\to\infty$.\qed

\subsection{The critical recovery rate}

In view of Proposition~\ref{P:infsurv}, we wish to find sufficient conditions
for $\prod_{k=0}^\infty(1-\xi(k))>0$. The next lemma casts the inductive
formula (\ref{induc2}) in a more tractable form.
\bl{\bf(Inductive formula)}\label{L:induc}
Let $\de(n),\al_k(n)$, and $\xi(n)$ be defined as in
Proposition~\ref{P:infsurv} and assume that the constants $(\al_k)_{k\geq 1}$
are positive. Set $\eps(k):=\de(k)/\al_1(k)$ $(k\geq 0)$. Then
$\xi(k)=f(1/\eps(k))$ and
\be
\eps(k+1)=\ffrac{\al_{k+1}}{\al_{k+2}}g(\eps(k))\qquad(k\geq 0),
\ee
where
\be\label{gdef}
g(\eps):=4\eps f(1/\eps)\quad(\eps>0),
\ee
and $f$ is the function defined in (\ref{fdef}).
\el
{\bf Proof} It is clear from (\ref{induc2}) that $\xi(n)=f(1/\eps(n))$ and
\be\label{alk}
\al_k(n)=2^{-n}\al_{k+n}\qquad(k\geq 1,\ n\geq 0).
\ee
Using (\ref{induc2}) once more, it follows that
\be
\eps(n+1)=\frac{2\xi(n)\de(n)}{\frac{1}{2}\al_2(n)}
=\frac{\al_1(n)}{\al_2(n)}\frac{\de(n)}{\al_1(n)}4f(\al_1(n)/\de(n))
=\frac{2^{-n}\al_{n+1}}{2^{-n}\al_{n+2}}
4\eps(n)f(1/\eps(n)).
\ee
\qed

\noi
The next lemma collects some elementary facts about the function $g$ from
Lemma~\ref{L:induc}.
\bl{\bf(The function $g$)}\label{L:g}
The function $g$ defined in (\ref{gdef}) is increasing on $(0,\infty)$ and
satisfies
\be\label{geps}
g(\eps)=8\eps^2+O(\eps^3)\qquad\mbox{as}\quad\eps\to 0.
\ee
\el
{\bf Proof} This follows from the fact that, by Lemma~\ref{L:gamma}, the
function $\eps\mapsto f(1/\eps)$ is increasing and satisfies
\be\label{feps}
f(1/\eps)=2\eps+O(\eps^2)\qquad\mbox{as}\quad\eps\to 0.
\ee
\qed

\noi
The next proposition answers the question when the infinite product in
(\ref{infbd}) is positive for $\de$ small enough.
\bp{\bf(Nontrivial survival bound)}\label{P:nontriv}
Let $(\al_k)_{k\geq 0}$ be nonnegative constants. For given $\de>0$, set
$\Pi(\de):=\prod_{k=0}^\infty(1-\xi(k))$, where the $(\xi(k))_{k\geq 0}$ are
defined as in Proposition~\ref{P:finsurv}. Then $\Pi(\de)$ is nonincreasing
in $\de$. Moreover, $\Pi(\de)>0$ for $\de$ suffiently small if and only if
\be\label{necsuf}
\sum_{k=m}^\infty2^{-k}\log(\al_k)>-\infty\quad\mbox{for some}\quad m\geq 0.
\ee
\ep
{\bf Proof} We start by showing that $\Pi(\de)$ is nonincreasing in $\de$. By
continuity, it suffices to prove this under the additional assumption that the
$\al_k$'s are all positive. In this case, we observe from Lemma~\ref{L:induc}
and the monotonicity of $g$ that the $\eps(k)$'s are nondecreasing in
$\de$. Since $\xi(k)=f(1/\eps(k))$ and $f$ is decreasing, it follows that the
$\xi(k)$'s are nondecreasing in $\de$, hence $\Pi(\de)$ is nonincreasing in
$\de$.

We next show that $\Pi(\de)>0$ for $\de>0$ suffiently small if and only if
(\ref{necsuf}) holds. If $\al_k=0$ for some $k\geq 1$, then
$\xi(k-1)=f(0)=\frac{1}{2}$, hence if infinitely many of the $\al_k$'s are
zero then $\Pi(\de)=0$ for all $\de>0$, while (\ref{necsuf}) is obviously
violated. If finitely many of the $\al_k$'s are zero, then we may start our
inductive formulas after the first $m$ iterations, where we observe that
$\de(m)$ can be made arbitrarily small by choosing $\de$ small enough. Thus,
without loss of generality, we may assume that the $\al_k$'s are all positive,
and under this assumption we need to show that $\Pi(\de)>0$ for $\de$
suffiently small if and only if
\be\label{necsuf2}
\sum_{k=0}^\infty2^{-k}\log(\al_k)>-\infty.
\ee
It is well-known that $\prod_{k=0}^\infty(1-\xi(k))>0$ if and only if
$\sum_{k=0}^\infty\xi(k)<\infty$. Using (\ref{feps}) and the fact that
$\xi(k)=f(1/\eps(k))$, it is easy to see that this is equivalent to
$\sum_{k=0}^\infty\eps(k)<\infty$.
\detail{Indeed, $\eps\leq f(1/\eps)\leq 3\eps$ for $\eps$ sufficiently small,
so $\sum_{k=0}^\infty\xi(k)<\infty$ if and only if
$\sum_{k=0}^\infty\eps(k)<\infty$.}

Now assume that (\ref{necsuf2}) holds, and, in view of (\ref{geps}),
define $(\ti\eps(k))_{k\geq 0}$ by
\be\label{tieps}
\ti\eps(0):=\frac{\de}{\al_1}\quad\mbox{and}\quad
\ti\eps(k+1):=9\ffrac{\al_{k+1}}{\al_{k+2}}(\ti\eps(k))^2\qquad(k\geq 0).
\ee
Then
\bc
\ti\eps(0)&=&\dis\ffrac{\de}{\al_1},\\
\ti\eps(1)&=&\dis9\ffrac{\al_1}{\al_2}(\ffrac{\de}{\al_1})^2,\\
\ti\eps(2)&=&\dis9\ffrac{\al_2}{\al_3}(9\ffrac{\al_1}{\al_2})^2
(\ffrac{\de}{\al_1})^4,\\
\ti\eps(3)&=&\dis9\ffrac{\al_3}{\al_4}(9\ffrac{\al_2}{\al_3})^2
(9\ffrac{\al_1}{\al_2})^4(\ffrac{\de}{\al_1})^8
=\frac{\frac{1}{9}(9\de)^8}{\al_4\al_3\al_2^2\al_1^4}.
\ec
More generally, it is not hard to see that
\be
\ti\eps(n)=\frac{\frac{1}{9}(9\de)^{2^n}}{\al_{n+1}
\prod_{k=1}^n(\al_k)^{2^{n-k}}}
\ee
By Lemma~\ref{L:suma} below, we can choose $\de$ sufficiently small such that
$\sum_{n=0}^\infty\ti\eps(n)<\infty$. By (\ref{geps}), there exists a $c>0$
such that $g(\eps)\leq 9\eps^2$ for all $\eps\leq c$. By making $\de$ smaller
if necessary, we can arrange that $\ti\eps(n)\leq c$ for all $n$, hence
$\sum_{n=0}^\infty\eps(n)\leq\sum_{n=0}^\infty\ti\eps(n)<\infty$.

On the other hand, assume that $\sum_{n=0}^\infty\eps(n)<\infty$ for some
$\de>0$ while (\ref{necsuf2}) does not hold. Define $(\ti\eps(k))_{k\geq 0}$
as in (\ref{tieps}) but with the factor $9$ replaced by $7$. By (\ref{geps}),
there exists a $c>0$ such that $7\eps^2\leq g(\eps)$ for all $\eps\leq
c$. Making $\de$ smaller if necessary, we can arrange that $\eps(n)\leq c$ for
all $n$, hence $\eps(n)\geq\ti\eps(n)$ for all $n$. Since
$\ti\eps(n)\to\infty$ by Lemma~\ref{L:suma} below, this leads to a
contradiction.\qed

\bl{\bf(Summability)}\label{L:suma}
For $\eta>0$, set
\be
F_\eta(n):=\frac{\eta^{2^n}}{\al_{n+1}
\prod_{k=1}^n(\al_k)^{2^{n-k}}}\qquad(n\geq 0).
\ee
If (\ref{necsuf2}) holds, then $\sum_{n=0}^\infty F_\eta(n)<\infty$ for
$\eta$ sufficiently small. On the
other hand, if (\ref{necsuf2}) does not hold, then
$\lim_{n\to\infty}F_\eta(n)=\infty$ for all $\eta>0$.
\el
{\bf Proof} We start by observing that
\be\ba{l}\label{limsup}
\dis\exists m\mbox{ s.t.\ }\forall n\geq m:\ F_\eta(n)<1\\[5pt]
\dis\quad\desd\quad\exists m\mbox{ s.t.\ }\forall n\geq m:\ 
\Big(2^n\log(\eta)-\log(\al_{n+1})
-\sum_{k=1}^n2^{n-k}\log(\al_k)\Big)<0\\[5pt]
\dis\quad\desd\quad\exists m\mbox{ s.t.\ }\forall n\geq m:\ \Big(\log(\eta)
-2^{-n}\log(\al_{n+1})-\sum_{k=1}^n2^{-k}\log(\al_k)\Big)<0\\[5pt]
\dis\quad\desd\quad\log(\eta)-\sum_{k=1}^\infty2^{-k}\log(\al_k)<0,
\ec
which is satisfied for $\eta$ sufficiently small if (\ref{necsuf2}) holds. In
this case, we may choose $\eta>0$ such that $K:=\sup_{n\geq
  0}F_\eta(n)<\infty$ and observe that for any $\eta'<\eta$
\be
\sum_{n=0}^\infty F_\eta(n)\leq
K\sum_{n=0}^\infty\big(\frac{\eta'}{\eta}\big)^{2^n}<\infty.
\ee
On the other hand, if (\ref{necsuf2}) does not hold, then a calculation as in
(\ref{limsup}) shows that for all $\eta>0$ there exists an $m$ such that for
all $n\geq m$ one has $F_\eta(n)>1$, and therefore, for any $0<\eta<\eta'$,
\be
\liminf_{n\to\infty}F_{\eta'}(n)
\geq\liminf_{n\to\infty}\big(\frac{\eta'}{\eta}\big)^{2^n}F_\eta(n)=\infty.
\ee
\qed

\subsection{Comparison argument}

{\bf Proof of Theorem~\ref{T:main}~(b)} For $N=2$,
Theorem~\ref{T:main}~(b) follows from Propositions~\ref{P:infsurv} and
\ref{P:nontriv}. To generalize this to arbitrary $N\geq 2$, we will use a
comparison argument.

Let $N\geq 2$ and let $X$ be a contact process on $\om_N$ with infection
rates as in (\ref{adef}) satisfying
\be\label{sumko}
\sum_{k=k_0}^\infty(N')^{-k}\log(\al_k)>-\infty\quad\mbox{for some }k_0\geq 1
\ee
where $N'=N$ in case $N$ is a power of two and $1<N'<N$ otherwise. For
notational convenience, we set $\ga_k:=\al_kN^{-k}$ $(k\geq 1)$, i.e.,
we let $(\ga_k)_{k\geq 1}$ denote the constants such that the infection rates
of $X$ are given by (compare (\ref{adef}))
\be\label{gadef}
a(i,j)=\ga_{|i-j|}\qquad(i,j\in\om_N,\ i\neq j).
\ee
Then, by (\ref{sumko}),
\be\label{gasum}
\sum_{k=k_0}^\infty(N')^{-k}\log(\ga_k)
=\sum_{k=k_0}^\infty(N')^{-k}\log(\al_k)-\sum_{k=k_0}^\infty(N')^{-k}k\log(N)
>-\infty.
\ee
We claim that we can choose $n,m\geq 1$ such that
\be\label{N2N}
(N')^m\leq 2^n\leq N^m.
\ee
If $N'=N$ is a power of two, then this is obviously satisfied with $m=1$ for
some $n\geq 1$. Otherwise, we observe that (\ref{N2N}) is equivalent to
\be
\frac{\log(N')}{\log 2}\leq\frac{n}{m}\leq\frac{\log(N)}{\log 2},
\ee
which is satisfied for some $m,n\geq 1$ since $N'<N$.

By an obvious monotone coupling, we can estimate $X$ from below by a contact
process $X'$ on $\om_N$ with infection rates of the form
$a'(i,j)=\ga'_{|i-j|}$, where
\be\ba{l}
\dis\ga'_1=\ga'_2=\cdots=\ga'_m=\min\{\ga_1,\ldots,\ga_m\},\\[5pt]
\dis\ga'_{m+1}=\ga'_{m+2}=\cdots=\ga'_{2m}
=\min\{\ga_{m+1},\ldots,\ga_{2m}\},\\[5pt]
\mbox{etcetera.}
\ec
Next, we consider a contact process $X''$ on $\om_2$ with infection
rates of the form $a''(i,j)=\ga''_{|i-j|}$, where
\be\ba{l}
\dis\ga''_1=\ga''_2=\cdots=\ga''_n=\ga'_m,\\[5pt]
\dis\ga''_{n+1}=\ga''_{n+2}=\cdots=\ga''_{2n}=\ga'_{2m},\\[5pt]
\mbox{etcetera.}
\ec
We claim that if $X''$ survives for a certain value of the recovery rate, then
so does $X'$. To see this, note that we can in a natural way identify $X'$
with a contact process on $\om_{N^m}$ with infection rates
$a(i,j)=\ga'_{m|i-j|}$. Likewise, we can in a natural way identify $X''$ with
a contact process on $\om_{2^n}$. Since $2^n\leq N^m$, we may regard
$\om_{2^n}$ as a subset of $\om_{N^m}$. Therefore, by surpressing infections
that go outside $\om_{2^n}$ we may estimate $X'$ from below by $X''$.

For $l\geq 0$, choose $i_l\in\{lm+1,\ldots,lm+m\}$ such that
\be
\ga_{i_l}=\min\{\ga_{lm+1},\ldots,\ga_{lm+m}\}.
\ee
Then, for $l\geq 0$ and $r=1,\ldots,n$, one has $\ga''_{ln+r}=\ga_{i_l}$,
hence by (\ref{N2N}) and (\ref{gasum}),
\be\ba{l}\label{gacc}
\dis\sum_{k=nl_0+1}^\infty2^{-k}|\log(\ga''_k)|
=\sum_{l=l_0}^\infty\sum_{r=1}^n2^{-(ln+r)}|\log(\ga''_{ln+r})|
=\sum_{l=l_0}^\infty\sum_{r=1}^n2^{-(ln+r)}|\log(\ga_{i_l})|\\[5pt]
\dis\quad\leq\sum_{l=l_0}^\infty2^{-ln}|\log(\ga_{i_l})|
\leq\sum_{l=l_0}^\infty(N')^{-lm}|\log(\ga_{i_l})|
\leq(N')^m\sum_{l=l_0}^\infty(N')^{-i_l}|\log(\ga_{i_l})|\\[5pt]
\dis\quad\leq(N')^m\sum_{i=ml_0+1}^\infty(N')^{-i}|\log(\ga_i)|<\infty
\ec
for some $l_0\geq 0$. If we write the infection rates of $X''$ in the form
$a''(i,j)=\al''_{|i-j|}N^{-|i-j|}$, then by (\ref{gacc}) and the calculation in
(\ref{gasum}) one has $\sum_{k=k_0}^\infty2^{-k}\log(\al''_k)>-\infty$ for
some $k_0\geq 0$, hence applying what we have already proved for $N=2$ we
conclude that $X''$ has a positive critical recovery rate and the same must be
true for $X$.  \qed

\appendix

\section{Coordinate reduction}\label{A:coord}

In this appendix we prove that Lemmas~\ref{L:onelev} and
\ref{L:twolev} imply formulas (\ref{RIi})~(i) and (ii),
respectively. The main problem is to invent good notation. Recall
that $S_n=\{0,1\}^{\om^n}$. For any $x\in S_n$ and $\De\sub\om^n$,
we let
\be
x\big|_\De:=(x(i))_{i\in\De}
\ee
denote the restriction of $x$ to $\De$. If $\De,\De'$ are disjoint
sets, $x\in\{0,1\}^\De$ and $x'\in\{0,1\}^{\De'}$, then we define $x\comb
x'\in\{0,1\}^{\De\cup\De'}$ by
\be
(x\comb x')(i):=\left\{\ba{ll}x(i)\quad&\mbox{if }i\in\De,\\
x'(i)\quad&\mbox{if }i\in\De'\ea\right.
\ee
For each $i\in\om^{n-1}$, we define $B_i\sub\om^n$ by (recall (\ref{block}))
\be
B(i):=B_1(i)=\{i'\circ i:i'\in\om^1\}.
\ee
Let $R,R'$ be as in Lemma~\ref{L:onelev}. Then we can write
\bc
\dis Rf(x)&=&\dis\sum_{x'\in S_1}R(x,x')f(x'),\\[5pt]
\dis R'f(y)&=&\dis\sum_{y'\in S_0}R'(y,y')f(y'),
\ec
where $R(x,x')$ and $R'(y,y')$ are the matrices of $R$ and $R'$, respectively.
We observe that
\bc
\dis R_if(x)&=&\dis\sum_{z\in\{0,1\}^{B(i)}}R(x|_{B(i)},z)
f(x|_{\om^n\beh B(i)}\comb z),\\[5pt]
\dis R'_if(y)&=&\dis\sum_{z\in\{0,1\}^{\{i\}}}R'(y(i),z)
f(y|_{\om^n\beh\{i\}}\comb z),
\ec
where we identify $\{0,1\}^{B(i)}\cong\{0,1\}^{\om^1}=S_1$ and
$\{0,1\}^{\{i\}}\cong\{0,1\}^{\om^0}=S_0$. Moreover,
\be
Pf(x)=\sum_{y\in\{0,1\}^{\om^{n-1}}}
\Big(\prod_{j\in\om^{n-1}}p(x_j,y(j))\Big)f(y).
\ee
Using the identification $x|_{B(i)}\cong x_i$, we calculate
\be\ba{l}
\dis R_iPf(x)=\sum_{z\in\{0,1\}^{B(i)}}R(x_i,z)
Pf(x|_{\om^n\beh B(i)}\comb z)\\[5pt]
\dis\quad=\sum_{z\in\{0,1\}^{B(i)}}R(x_i,z)
\sum_{y\in\{0,1\}^{\om^{n-1}}}
\Big(\prod_{j\in\om^{n-1}}p((x|_{\om^n\beh B(i)}\comb z)_j,y(j))\Big)
f(y)\\[5pt]
\dis\quad=\sum_{y\in\{0,1\}^{\om^{n-1}}}
\sum_{z\in\{0,1\}^{B(i)}}R(x_i,z)p(z,y(i))
\Big(\prod_{j\in\om^{n-1}\beh\{i\}}p(x_j,y(j))\Big)f(y)\\[5pt]
\dis\quad=\sum_{y\in\{0,1\}^{\om^{n-1}}}
\sum_{z\in\{0,1\}^{\{i\}}}p(x_i,z)R'(z,y(i))
\Big(\prod_{j\in\om^{n-1}\beh\{i\}}p(x_j,y(j))\Big)f(y)\\[5pt]
\dis\quad=\sum_{y\in\{0,1\}^{\om^{n-1}}}
\sum_{z\in\{0,1\}^{\{i\}}}p(x_i,y(i))R'(y(i),z)
\Big(\prod_{j\in\om^{n-1}\beh\{i\}}p(x_j,y(j))\Big)
f(y|_{\om^{n-1}\beh\{i\}}\comb z)\\[5pt]
\dis\quad=\sum_{y\in\{0,1\}^{\om^{n-1}}}
\Big(\prod_{j\in\om^{n-1}}p(x_j,y(j))\Big)
\sum_{z\in\{0,1\}^{\{i\}}}R'(y(i),z)f(y|_{\om^{n-1}\beh\{i\}}\comb z)=PR'_if.
\ec
Here we have used Lemma~\ref{L:onelev} in the fourth equality. In the
fifth equality, we have reordered our sums by relabelling $y(i)$ and
$z$.

The formal proof of formula (\ref{RIi})~(ii) is similar, but even more
cumbersome. Letting $I$ and $I'_x$ be as in Lemma~\ref{L:twolev}, we
can write, in matrix notation,
\bc
\dis If(x)&=&\dis\sum_{x'\in S_2}I(x;x')f(x')
=\sum_{z\in S_1}\sum_{z'\in S_1}I(x_0,x_1;z,z')f(z,z'),\\[5pt]
\dis I'_xf(y)&=&\dis\sum_{y'\in S_1}I'_x(y;y')f(y')
=\sum_{z\in\{0,1\}}\sum_{z'\in\{0,1\}}I_{x_0,x_1}(y(0),y(1);z,z')f(z,z').
\ec
Then
\bc
\dis I_{ij}f(x)&=&\dis\sum_{z\in\{0,1\}^{B(i)}}\sum_{z'\in\{0,1\}^{B(j)}}
I(x|_{B(i)},x|_{B(j)};z,z')
f(x|_{\om^n\beh(B(i)\cup B(j))}\comb z\comb z'),\\[5pt]
\dis\ov I_{ij}f(x,y)&=&\dis
\sum_{z\in\{0,1\}^{\{i\}}}\sum_{z'\in\{0,1\}^{\{j\}}}
I'_{x|_{B(i)},x|_{B(j)}}(y(i),y(j);z,z')
f(y|_{\om^{n-1}\beh\{i,j\}}\comb z\comb z'),
\ec
and
\be
\ov Pf(x)=\sum_{y\in\{0,1\}^{\om^{n-1}}}
\Big(\prod_{k\in\om^{n-1}}p(x_k,y(k))\Big)f(x,y).
\ee
Using the fact that $x|_{B(i)}\cong x_i$, we calculate
\be\ba{l}
\dis I_{ij}Pf(x)=\sum_{z\in\{0,1\}^{B(i)}}\sum_{z'\in\{0,1\}^{B(j)}}
I(x_i,x_j;z,z')
Pf(x|_{\om^n\beh(B(i)\cup B(j))}\comb z\comb z')\\[5pt]
\dis\quad=\sum_{z\in\{0,1\}^{B(i)}}\sum_{z'\in\{0,1\}^{B(j)}}
I(x_i,x_j;z,z')\\
\dis\qquad\quad\cdot\sum_{y\in\{0,1\}^{\om^{n-1}}}\Big(\prod_{k\in\om^{n-1}}
p((x|_{\om^n\beh(B(i)\cup B(j))}\comb z\comb z')_k,y(k))\Big)f(y)\\[5pt]
\dis\quad=\sum_{y\in\{0,1\}^{\om^{n-1}}}
\sum_{z\in\{0,1\}^{B(i)}}\sum_{z'\in\{0,1\}^{B(j)}}
I(x_i,x_j;z,z')p(z_i,y(i))p(z_j,y(j))\\
\dis\qquad\quad\cdot
\Big(\prod_{k\in\om^{n-1}\beh\{i,j\}}p(x_k,y(k))\Big)f(y)\\[5pt]
\dis\quad=\sum_{y\in\{0,1\}^{\om^{n-1}}}
\sum_{z\in\{0,1\}^{\{i\}}}\sum_{z'\in\{0,1\}^{\{j\}}}
p(x_i,z)p(x_j,z')I'_{x_i,x_j}(z,z';y(i),y(j))\\
\dis\qquad\quad\cdot
\Big(\prod_{k\in\om^{n-1}\beh\{i,j\}}p(x_k,y(k))\Big)f(y)\\[5pt]
\dis\quad=\sum_{y\in\{0,1\}^{\om^{n-1}}}
\sum_{z\in\{0,1\}^{\{i\}}}\sum_{z'\in\{0,1\}^{\{j\}}}
p(x_i,y(i))p(x_j,y(j))I'_{x_i,x_j}(y(i),y(j);z,z')\\
\dis\qquad\quad\cdot
\Big(\prod_{k\in\om^{n-1}\beh\{i,j\}}p(x_k,y(k))\Big)
f(y|_{\om^{n-1}\beh\{i,j\}}\comb z\comb z')\\[5pt]
\dis\quad=\sum_{y\in\{0,1\}^{\om^{n-1}}}
\Big(\prod_{k\in\om^{n-1}}p(x_k,y(k))\Big)\\
\dis\qquad\quad\cdot
\sum_{z\in\{0,1\}^{\{i\}}}\sum_{z'\in\{0,1\}^{\{j\}}}
I'_{x_i,x_j}(y(i),y(j);z,z')
f(y|_{\om^{n-1}\beh\{i,j\}}\comb z\comb z')\\[5pt]
\dis\quad=\ov P\,\ov I_{ij}f(x).
\ec
Here we have used Lemma~\ref{L:twolev} in the fourth equality, and in
the fifth equality, we have reordered our sums by relabelling
$y(i),y(j),z$, and $z'$.

\subsubsection*{Acknowledgements}

We thank the referees for a careful reading of the manuscript, and B\'alint
T\'oth, Deepak Dhar and Roman Kotecky for pointing out the references
\cite{Dys69} and \cite{BM87}.

\newcommand{\noopsort}[1]{}

\end{document}